\documentclass[12pt]{amsart}

\usepackage{amsmath,amssymb,amsthm,epsfig}
\usepackage{enumerate}
\usepackage{xcolor}
\usepackage{bbm}
\usepackage{appendix}
\usepackage[margin=20truemm]{geometry}
\usepackage[colorinlistoftodos]{todonotes}
\usepackage{version}
\usepackage[all]{xypic}

\newif\ifshowrevisions
\showrevisionstrue
\ifshowrevisions
  
  \newenvironment{revblock}{\begingroup\color{red}}{\endgroup}
\else

\fi

\DeclareFontFamily{OT1}{rsfs}{}
\DeclareFontShape{OT1}{rsfs}{n}{it}{<-> rsfs10}{}
\DeclareMathAlphabet{\curly}{OT1}{rsfs}{n}{it}

\newtheorem{Thm}{Theorem}[section]

\newtheorem{cor}[Thm]{Corollary}

\newtheorem{prop}[Thm]{Proposition}

\newtheorem{``Conj"}[Thm]{``Conjecture"}

\newtheorem{Ques}[Thm]{Question}
\newtheorem*{MQues}{Main Question}

\theoremstyle{remark}
\newtheorem{Rem}[Thm]{Remark}
\newtheorem{ex}[Thm]{Example}
\theoremstyle{definition}
\newtheorem{defn}[Thm]{Definition}

\newtheorem{Setup}{Setup}

\newcommand{\Spec}{\mathop{\mathrm{Spec}}\nolimits}

\def\dim{\operatorname{dim}}

\def\vol{\operatorname{vol}}

\def\R{\mathbb{R}}

\def\C{\mathbb{C}}
\def\G{\mathbb{G}}
\def\Z{\mathbb{Z}}
\def\O{\mathcal{O}}
\def\P{\mathbb{P}}
\def\Q{\mathbb{Q}}

\def\A{\mathbb{A}}

\def\X{\mathcal{X}}
\def\Y{\mathcal{Y}}

\def\D{\mathcal{D}}
\def\W{\mathcal{W}}

\def\M{\mathcal{M}}

\excludeversion{memo}

\begin{document}

\title[Algebraic geometry of bubbling]
{Algebraic geometry of bubbling K\"ahler metrics}
\author{Yuji Odaka}

\maketitle
\thispagestyle{empty}

\setcounter{tocdepth}{1} 
\noindent
\tableofcontents

\begin{abstract}
We give an algebro-geometric framework to study 
bubbling phenomena of K\"ahler metrics with Euclidean volume growth, 
after \cite{DSII, Song, dBS}. 
That is, for a degenerating family 
of log terminal (canonical) Gorenstein singularities 
we give an algorithm to construct a finite 
sequence of birational modifications of the family with milder 
degenerations, 
and partially relate to analytic bubbling 
constructions in \cite{Song, dBS}. 
We also 
provide more explicit equivalent approaches via coordinates. 
\end{abstract}

\section{Introduction}\label{sec:intro}

We make a purely algebro-geometric study of  
degenerations of normal singularities in this paper, 
motivated by the recent works on asymptotic behaviour of K\"ahler metrics \cite{DSII, Song, dBS}. 

For simplicity, we work over the field of 
complex numbers $\C$. 
Our main setup is as follows. 

\begin{Setup}\label{Setup}
Let $\pi\colon \mathcal{X}\twoheadrightarrow \Delta\ni 0$ be as follows: 
$\Delta$ is a smooth affine curve (germ) over $\mathbb{C}$ with a closed point $0$, 
$\pi$ is a flat affine 
morphism between $\mathbb{C}$-varieties 
of relative dimension $n$. We fix a uniformizer $t$ of $\mathcal{O}_{\Delta,0}$ and 
denote the base change $\mathcal{X}\times_\Delta {\rm Spec}\C[[t]]$ as $\X_{\C[[t]]}$. 
(Occasionally, we have to take ramified base change 
$\Delta^{(d)}\to \Delta$
of degree $d\in \Z_{>0}$, 
which lifts to ${\rm Spec}\mathbb{C}[[t^{1/d}]]\to {\rm Spec}\mathbb{C}[[t]]$.) 
For simplicity of the exposition, 
we also assume $\pi$ is a Gorenstein morphism so that 
$K_{\X}$ is Cartier and each fiber 
is log terminal (and Gorenstein again). 
\footnote{
These Gorenstein conditions should not be essential, just put to avoid 
complication of dealing with $\Q$-Gorenstein deformation 
space 
in the proof of Theorem \ref{bblcst:def} (cf., e.g., 
\cite[Chapter 9]{Kol23}). 
See Remark \ref{ll.rem}. Also, 
recall that normal Gorenstein singularity is 
log terminal iff canonical iff rational. }
We set $X:=\X_0:=\pi^{-1}(0)$. 
We further fix an algebraic section $\sigma\colon \Delta \to \X$  
with $\sigma(0)=x$. 
\end{Setup}

In section \ref{sec:2}, we algebraically construct 
(almost) canonical, finite ordered 
sequences of modifications of $\X$ with milder degenerations. 
For instance, if the general fibers are 
$2$-dimensional 
smooth, degenerating to 
$X$ with a $2$-dimensional $A_k$-singularity, 
then we obtain finite other fillings with 
$A_i$-singularities with smaller Milnor 
number $i$ each (cf., \cite{dBS}) and terminate at 
smooth filling. The process is 
essentially different from the classical 
simultaneous resolution of Brieskorn-Grothendieck type. 
Morally speaking, our algorithm gives 
canonical ``limit singularities" of various rescales beyond the given $X$ itself or its (metric) tangent cone, 
reflecting the whole formation process of the singularities of 
$X$, in the context of K\"ahler-Einstein geometry and the 
K-stability. 
We explain its 
detailed meaning later; 
here is a brief overview of the 
differential geometric background. 

Suppose we have a ``good" continuous 
family of (singular) K\"ahler metrics $g_t$ on $\X_t=\pi^{-1}(t)$ 
in the pluri-potential theoretic sense (cf., e.g., 
\cite{BBEGZ}), for $t\neq 0$, 
especially of K\"ahler-Einstein metrics, we would like to consider and answer 
the following question: 

\begin{MQues}[see \cite{Song, dBS}]
How can we determine the set of possible non-collapsing pointed Gromov-Hausdorff limits 
of $(x_t\in \X_t, c_t^2 g_t)$ for $t\to 0$ where $c_t>0$ is a sequence of diverging 
rescaling parameters? 
\end{MQues}

In this paper, 
following \cite{Song, dBS}, 
we call such a non-collapsing pointed Gromov-Hausdorff limit simply {\it bubbles} 
of 
the family $(\X_t,g_t)$ for $t\to 0$. 
By the groundbreaking work \cite[1.4]{DSII}, extending some classical cases, 
such limits often naturally admit the 
structure of an affine log terminal 
variety with (weak) Ricci-flat K\"ahler metric. 
Indeed, our aforementioned main construction of affine varieties 
is expected to underlie these bubbles. 
They have Euclidean volume growth simply by the 
Bishop's inequality. 
This is a (K\"ahler) metric version of so-called {\it bubbling} phenomenon since \cite{SU}, 
as later also studied by 
\cite{Nakajima, Anderson, Bando, BKN} in mainly $4$-dimensional Riemannian geometric setup, 
which in particular found out the finite energy (the $L^2$-norm of the curvature) rescaling provides smooth 
asymptotically locally Euclidean (ALE) manifolds 
as the deepest bubbles. 
Then, recently we witnessed a breakthrough \cite{DSII, 
Song, dBS} in the higher dimensional K\"ahler setup. 
In particular, \cite{Song} proved the existence of ``minimal" bubble 
in a differential geometric sense, finiteness of bubbles, 
and raised Conjecture 5.4 in {\it loc.cit} to expect a $2$-step 
algebraic construction of the underlying varieties of the minimal bubbles, 
and \cite{dBS} showed many examples and a refinement of the conjecture. 

Motivated in part by {\it loc.cit}, 
this paper aims to provide a systematic algebro-geometric or non-archimedean 
approach to the problem. 

The main contents of this paper provide an algebro-geometric
construction motivated by the conjectural picture
in \cite{Song,dBS}. 
The following is a brief outlook of our statements 
for our algorithm, but see the detailed statements 
in later referred subsections. 
\begin{Thm}\label{intro.thm}
Consider an arbitrary flat degeneration of pointed 
$n$-dimensional log terminal Gorenstein  
singularities 
$\pi\colon \X(\supset \sigma(\Delta))\to \Delta\ni 0$ 
(see Setup \ref{Setup}), 
which is strictly degenerating in a certain sense (see Theorem \ref{bblcst:def}). Then, the following hold. 
\begin{enumerate}
\item \label{intro.thm.1}
(Theorem \ref{bblcst:def}) There is a purely algebro-geometric 
$2$-step construction
\footnote{the following symbol $\rightsquigarrow$ 
means each step of the construction, and is not a map.}
of algebro-geometric version of 
minimal 
\footnote{The minimality is expected to correspond to that of 
the diverging order of the 
rescaling constants $\{c_t\}_t$ 
(cf., also \cite{Song} and Definition \ref{DG.min}), 
but can be characterized algebro-geometrically; in terms of 
the local normalized volumes 
(at the base points). We clarify the detailed meaning later.}
bubbles as affine varieties 
$$(\sigma(\Delta)\subset X\subset \X) \rightsquigarrow \X'_{\rm min,0}\rightsquigarrow 
\X''_{\rm min,0}.$$ 
Here, 
$\X'_{\rm min,0}, \X''_{\rm min,0}$ 
are both still $n$-dimensional affine varieties 
with only log terminal Gorenstein singularities. 

We call 
$\X'_{\rm min,0}$ (resp., $
\X''_{\rm min,0}$) 
{\it (algebraic) semistable minimal bubble}  
(resp., 
{\it (algebraic) polystable minimal bubble}) 
for simplicity. 
\footnote{\cite{Song} (resp., \cite{dBS}) conjectured 
such construction and called the semistable 
minimal bubble 
$\X'_{\rm min,0}$ 
under different names as 
algebraic minimal bubble (resp., 
(minimal) weighted bubble). 
Note that our construction of 
$\X''_{\rm min,0}$ is also algebraic.}

\item \label{intro.thm.4}
(Theorem \ref{bblcst:def}, Theorem \ref{kl}) 
The first step construction 
$$X\subset \X \rightsquigarrow \X'_{\rm min,0}$$ 
can be described in terms of weighted blow ups. 
Conditionally, this construction of 
$\X'_{\rm min,0}$ is canonical. 
\end{enumerate}
We refer to these obtained (non-canonical) total spaces 
of the new degenerations to 
$\X'_{\rm min,0}$ 
and  $\X''_{\rm min,0}$ respectively 
as 
$\X'_{\rm min}$ and $\X''_{\rm min}$. 
By replacing $\X$ by $\X''_{\rm min}$, and 
iteration of such procedure, we obtain further (deeper) modifications of $\X$. 
\end{Thm}

This whole 
process can be seen as a certain algorithm 
of improving the degenerations with some differential 
geometric motivation. Its double iterative nature 
may make it look too complicated or even artificial, 
but it is unavoidable 
and indeed parallels the actual bubbling structure 
and hence we believe it is the nature of the arts in general. 
Nevertheless, 
it becomes simpler in many classical examples. 
For instance, in the $1$ or $2$-dimensional 
classical cases considered in 
\cite[\S 2.3, 3.3, 4.2]{dBS}, 
the second step of \eqref{intro.thm.1} each time 
is trivial so that 
$\X'_{\rm min,0}=\X''_{\rm min,0}$ recovers the 
differential geometric construction i.e., 
a certain pointed Gromov-Hausdorff limit 
(at the underlying affine varieties level) and the 
above construction of $\X'_{\rm min,0}$ 
is canonical. 
Thus, by iterations, we obtain the whole set of 
possible differential 
geometric bubbles along $\X\supset \sigma(\Delta)$ 
with the local metrics given in {\it loc.cit}. 
See Theorem \ref{AGDG.cl}. 

We leave future work the precise relation between our 
construction and the analytic construction of \cite{Song} to future work; below we retain 
direct comparisons only in the classical cases treated in \cite{dBS}. 

Although our construction is fairly general and may look abstract at first, 
it is based on 
various explicit calculations as experiments. Indeed, our section 
\ref{sec:min1} connects to weighted blow ups of the total space of deformations. 
We review one of the simplest 
examples after \cite{dBS} here. 
In this case, we only need a single process rather than 
the $2$-step processes i.e., $\X'_{\rm min}=\X''_{\rm min}$, 
and they are achieved only after a finite base change of $\Delta$. 

\begin{ex}[cf., {\cite{BKN, Bando}, \cite[\S 3.3]{dBS} and Theorem \ref{AGDG.cl}}]\label{ex.dBS}
Consider a degeneration of pointed affine smooth surfaces to a 
$A_3$-singular surface as 
$\X\subset \Delta_t\times \C^3_{x,y,z}$ cut out by the equation 
\begin{align}
xy=z(z-t)(z-t^2)(z-t^3)
\end{align}
whose base points are given by the $0$-section $\sigma$. 
Here, the subindices $x,y,z,t$ are 
coordinates. 

We can consider a family of hyperK\"ahler metrics for  
each $t$ by the Gibbons-Hawking ansatz with the base 
$\R^3$ and monopole points $x_1, x_2, x_3, x_4$ and 
harmonic functions of the form $\frac{1}{2}\sum_i \frac{1}{|x-
x_i|}$. 
In this case, 
\cite[\S 3.3]{dBS} 
describes how it bubbles at $t\to 0$ around the origin, 
and coincides with the algebraic procedure 
Theorem \ref{intro.thm}, and is explicit as follows. 
(Indeed, our general construction 
is originally obtained as a generalization of their 
examples). 

Firstly, we consider the weighted blow up 
of $(0,0,0,0)\in \X$ 
with weights $(2,2,1,1)$ and removal of the 
strict transform of the central fiber $\X_0$, which  
gives $\X'_{\rm min}=\X''_{\rm min}$ with the equation 
\begin{align}
xy=z(z-1)(z-t)(z-t^2), 
\end{align}
after replacing $\frac{x}{t^2}$ as $x$, 
$\frac{y}{t^2}$ as $y$, $\frac{z}{t}$ as $z$ while keeping $t$. 
This is a degeneration to 
$xy=z^3(z-1)$ when $t\to 0$, with 
$A_2$-singularity at 
$(0,0,0,0)\in \Delta_t\times \C_{x,y,z}^3$. 
This admits a ALE hyperK\"ahler orbifold metric, 
again obtained via the Gibbons-Hawking ansatz over $\R^3$ 
(cf., \cite[\S 3.3]{dBS}). 

Then after the base change with respect to 
the double covering 
$\Delta_{\sqrt{t}}\to \Delta_t$, 
we replace $x, y,z,t$ by 
$\frac{x}{t^{1.5}}, \frac{y}{t^{1.5}},\frac{z}{t}, t^{0.5}$ 
and still denote them by $x,y,z,t$ for simplicity of notation. 
Then, the weighted blow up of the base changed family with 
weights $(3,3,2,2)$ (again with the strict transform of the central fiber removed) is defined by 
\begin{align}
xy&=z(t^2 z-1)(z-1)(z-t^2)\\
  &=(t^2 z-1)(z-1)\cdot z(z-t^2), 
\end{align}
which is a degeneration to a 
ALE hyperK\"ahler orbifold bubbling with a 
$A_1$-singularity where the metric is infinitesimally 
$\C^2/\pm$. 

As the final step, 
we again do a weighted blow up 
with weights $(2,2,2,1)$ for the coordinates $(x,y,z,t)$ 
and remove the strict transform of the central fiber. Then, 
by replacing the coordinates by new 
$(\frac{x}{t^2},\frac{y}{t^2},\frac{z}{t^2},t)$ 
which we denote again by $(x,y,z,t)$, 
we obtain the new equation 
\begin{align}
    xy=(t^4 z-1)(t^2z-1)\cdot z(z-1), 
\end{align}
inside $\Delta_t\times \C^3_{x,y,z}$. 
Note that this is finally a smooth family and we have an 
asymptotically conical 
hyperK\"ahler metric of Eguchi-Hanson type at the central fiber. 
\end{ex}

As we clarify, in general, we need more complicated processes 
due to possible irregularity of the corresponding Sasaki-Einstein 
link i.e., irrationality of its Reeb vector field. 

We make this kind of explicit study as Example \ref{ex.dBS} 
via coordinates more generally. 
That is, after the proof of Theorem \ref{intro.thm} ($=$Theorems \ref{bblcst:def}), 
we propose 
more explicit ways of constructions of the candidates of bubbles, 
via coordinates/equations (resp., 
valuations) in \S \ref{sec:min1} 
(resp., \cite{Odfuturebble}). 
We conjecture that they coincide with 
the construction of Theorem \ref{bblcst:def}. 
\S \ref{sec:nwd} 
may be of own independent 
interest, 
which are respectively 
the systematic treatment of negative weight deformation 
originated in \cite{CH, SZ}. 

From some personal perspectives, this paper 
is also indirectly motivated by hope to deepen our understanding of 
the projective (resp., quasi-projective) K-moduli of Fano varieties 
(resp., polarized Calabi-Yau varieties) in the sense of \cite[5.2]{Od10}, after recent completion of the 
construction (cf., e.g., \cite{BlumLiu, BHLLX, LXZ} resp., \cite{Bir, Od21}). 
Recall that those K-moduli were constructed and understood before {\it loc.cit} equivalently as 
Gromov-Hausdorff (partial) compactifications 
from more differential geometric viewpoints 
(cf., e.g., \cite{OSS, SSY, OdF, LWXF, Bir, Od21}), 
following the breakthrough work of \cite{DS}. 
Thus it is natural to hope this work on 
bubblings will give a deepening of our 
understanding of these moduli spaces. 

\section{Minimal bubble}\label{sec:2}

\subsection{Review of the basics}\label{sec:rev}

We work in Setup \ref{Setup} in the introduction. 
In this section, we 
give algebraic approaches to the ``minimal bubble" in the 
sense of \cite{Song, dBS}, whereas the minimality 
roughly refers to that of the order of rescaling parameters $c_t$ for $t\to 0$. 
We review the definition for convenience. 

\begin{defn}[{\cite[Definition 1.9]{Song}}]\label{DG.min}
In the context of Setup \ref{Setup}, 
its (differential geometric) {\it minimal 
bubble} refers to the pointed Gromov-Hausdorff limits of 
$(\sigma(t)\in \X_t, c_t^2 g_t)$ for certain diverging sequence $c_t>0$, 
which is not a Riemannian metric cone, but whose tangent cone at infinity (\cite[\S 3.4]{DSII}) is 
the metric local tangent cone at $\sigma(0)\in X$, which we denote as $C_{\sigma(0)}(X)$. 
\end{defn}

Firstly, we discuss several a priori distinct 
algebraic approaches to 
its constructions, in different subsections, and will discuss the conjectural equivalence 
later. We crucially use the theory of local normalized volume of 
log terminal singularities due to 
\cite{Li, Liu} (cf., also \cite{Od}), with its origin in \cite{DSII}. 
We briefly review the definition and properties we use here. 
Its smooth case can be traced back to much 
earlier \cite{Mustata, dFEM} (and also related \cite{ELS}) and 
more detailed reviews are available at \cite{LLX}, \cite[\S 2]{Od24a} 
for instance.

\begin{defn}[\cite{ELS}]
For any quasi-monomial valuation $v$ of 
$n$-dimensional 
log terminal variety 
$X$ centered at a closed point $x\in X$, 
we define the volume ${\rm vol}_x(v)$ as 
$\limsup_{m\to \infty} \dfrac{\dim_k(\O_{X,x}/\{f\mid v(f)\ge m\})}{m^n/n!}.$
This $\limsup$ is known to be $\lim$ (cf., {\it loc.cit}). 
\end{defn}

Then, the normalized volume $\widehat{{\rm vol}}_x(v)$
is defined as $A_X(v)^n {\rm vol}_x(v)$ by \cite{Li}, 
where $A_X(v)$ is the log discrepancy in the generalized sense by \cite{JM}. 
Then the following is known by now. In this paper, we do not make any substantial 
contribution but rather use it. 
See more references in 
\cite{LLX}, \cite[\S 2]{Od24a}. 
Morally speaking, this process is to modify a ``locally K-stable" i.e., log terminal singularity 
(cf., \cite{Od, GT}) 
to cone, which is still (K-semi/poly)stable. 

\begin{Thm}[Scale up conification cf., \cite{DSII, Li, Liu, Blum, HS, LX, LWX, LLX, XZ2, Od24a}]\label{lconif}
For any log terminal singularity $x\in X$, 
there is a unique quasi-monomial 
valuation $v=v_X$ with the center $x$ 
up to the $\R_{>0}$-action, which minimizes 
the normalized volume $\widehat{\vol}(-)$. 
\begin{enumerate}
\item Then, the valuation 
$v_X$ coincides with the vanishing order function in \cite[\S 3.2]{DSII} 
for the local K\"ahler metric in the setup of {\it loc.cit}. 
We denote the rational rank of $v_X$ as $r$. 
Further, $\widehat{\rm vol}(v_X)/n^n$ has a differential geometric meaning as 
the volume density of the link of the local metric tangent cone. 
We denote the minimized value 
$\widehat{\rm vol}(v_X)$ 
as $\widehat{\rm vol}(x\in X)$ or 
$\widehat{\rm vol}_X(x)$, in this paper. 

\item 
From a more algebraic viewpoint, the graded ring 
$${\rm gr}_v(\mathcal{O}_{X,x})
:=\oplus_{s\in {\rm Im}(v_X)}
(\{f\in \mathcal{O}_{X,x}\mid v(f)\ge s\}/
\{f\in \mathcal{O}_{X,x}\mid v(f)> s\})$$ 
is a finite type $\mathbb{C}$-algebra, and 
$$W:=\Spec({\rm gr}_v(\mathcal{O}_{X,x}))$$
is a K-semistable Fano cone equipped with the natural action of the algebraic torus 
$T\simeq \G_m^r$. 
\item \label{rev:pwd}
For a certain rational polyhedral cone $\tau$ of the dual of 
the groupification of ${\rm Im}(v_X)$, there is a positive weight 
degeneration of $X$ to $W$, over $U_\tau$ which is a 
$\Q$-Gorenstein family. 

\item There is a (non-canonical) affine 
test configuration of 
$W$ to a unique K-polystable Fano cone $C$ with a (weak) 
Ricci-flat K\"ahler cone metric. In the setup of \cite[\S 2-\S 3.3]{DSII}, 
it recovers the local metric tangent cone of local K\"ahler metric 
(as conjectured in {\it loc.cit}). 

\item (\cite[Proposition 4.1]{NO}) $W$ and $C$ have the same $\Q$-Gorenstein index as that of 
$x\in X$. 
\end{enumerate}
\end{Thm}

In particular, if $x\in X$ is Gorenstein, so are $W$ and $C$. 


\subsection{Construction as deformation of cone}\label{sec:min0}

We use 
theorems in \cite{Od24b} 
after \cite{Langton, HL, AHLH} to 
algebraic construction of the conjectural (minimal) bubble 
and discuss it. We also discuss 
a reformulation of 
negative weight deformations 
in the sense of \cite{CH}, 
after \cite{Od24b}, in which we also called it 
{\it scale down deformation}. 

\subsubsection{Setup and algebraic construction}

In the following two subsubsections, 
we set up stages and list some statements in more precise form  
we prove in 
the later subsubsections. 

\begin{Setup}\label{Setup2}
In the context of Setup \ref{Setup}, 
as Example 2.8 in \cite{Od24b}, denote by 
$v_X$ the valuation of $\O_{X,x}$ whose center is $x$ and 
minimizes the normalized volume $\widehat{\vol}(-)$ (\cite{DSII, Li, Blum, XZ2}). 

We set the groupification of ${\rm Im}(v_X)$ as $M\simeq \Z^{r}$, 
$N:=M^*:={\rm Hom}_{\Z}(M,\Z)$ denotes the 
dual lattice, 
and set 
$T:=N\otimes \G_m=Hom(M,\G_m)\simeq \G_m^r$, 
where $r$ is the $\Q$-rank of $v_X$. 
\end{Setup}



The following main theorem \ref{bblcst:def} and its coordinate 
description 
aim at algebraic reconstruction and generalization of 
the bubbles via metrics in \cite{DSII, Song}. 
By using the general framework of \cite{Od24b}, we can make it more 
canonical and general by working over more 
general bases 
(cf., earlier version on arXiv if interested), 
but we provide a simpler form to avoid the complication. 

\begin{Thm}[Algebro-geometric construction of (minimal) bubbles]\label{bblcst:def}
Under the Setup \ref{Setup} $\land$ \ref{Setup2},  
we further assume that 
$\sigma(t)\in X_t$ strictly degenerates in the sense that 
$$\widehat{\vol}_{\X_0}(\sigma(0))<\widehat{\vol}_{\X_t}(\sigma(t))$$ for any 
$t\neq 0$. 
Here $\widehat{\vol}(-)$ denotes the normalized volume. 

Then, possibly after finite base change of $\pi$, 
there is a pair of 
birational transforms 
$\X\dashrightarrow \X_{\rm min}'\to \Delta$ and 
$\X\dashrightarrow \X_{\rm min}''\to \Delta$ as 
affine faithfully flat morphisms, (relatively/fiberwise) log terminal 
Gorenstein families, which satisfy the following. 

We denote the projection $\X_{\rm min}'\to \Delta$ as $\pi'$, 
and $\X_{\rm min}''\to \Delta$ as $\pi''$ respectively. 

\begin{enumerate}
\item $\sigma|_{\Delta\setminus \{0\}}$ extends to a holomorphic section 
$\Delta \to \X_{\rm min}'$, denoted by $\sigma'$ (resp., to a holomorphic section 
$\Delta \to \X_{\rm min}''$, denoted by $\sigma''$). 
\item \label{vol.ineq}
For the central fibers $\X'_{\rm min,0}$ of $\X$ (resp., $\X''_{\rm min,0}$ of $\X''_{\rm min}$), 
we have 
\begin{align}\tag{*}\label{vol.ineq2}
 \widehat{\vol}_{\X'_{\rm min,0}}(\sigma'(0))= 
\widehat{\vol}_{\X''_{\rm min,0}}(\sigma''(0))>
\widehat{\vol}_{\X_0}(\sigma(0)). 
\end{align}

\item \label{Kss.Song}
The central fiber $\X'_{{\rm min},0}:=\pi'^{-1}(0)$ of 
$\X_{\rm min}'\to \Delta$, 
is an algebraic 
negative weight deformation (in the sense of \cite[\S 2]{Od24b}, after \cite[Definition 1]{CH})
of a K-semistable $\Q$-Fano cone 
$(W,\xi)$ over $\mathbb{A}_{\mathbb{C}}^1$ 
in a canonical manner, 
which we denote as 
$\pi_\mathcal{W}\colon \mathcal{W}\twoheadrightarrow \A_{\mathbb{C}}^1$. 
In particular, 
$\pi_\W^{-1}(0)=W$ and 
$\pi_\W^{-1}(1)\simeq \X'_{{\rm min}, 0}$. 

\item \label{Kss.Song2}
$W=\Spec(gr_{v_X}(\O_{X,x}))$ with the natural 
action of the algebraic torus 
$T=\G_m \otimes_{\Z}N$ 
where $N$ denotes the dual of the groupification of the 
valuation semigroup of the volume-minimizing valuation $v_X$. 
That is, $v_X$ is the valuation of $\O_{X,x}$ whose center is $x$ and 
minimizes the normalized volume $\widehat{\vol}(-)$ 
and 
$gr_{v_X}(\O_{X,x})$ is the associated graded ring.

\item \label{Kps.Song}
The central fiber 
$\X''_{{\rm min},0}:=(\pi'')^{-1}(0)$ of $\X''_{\rm min}$ 
is an algebraic negative weight deformation 
\footnote{The notion is introduced in 
\cite[\S 2]{Od24b}, after preceding notion by 
\cite[1.7]{CH}. We clarify their relation in 
\S \ref{sec:nwd}. 
There is also another preceding related discussion in 
\cite[\S 6]{SunZhang}.). }\label{ft5}
(\cite[\S 2]{Od24b}, after \cite{CH}) 
of a K-polystable $\Q$-Fano cone $(C,\xi_C)$ in the sense of \cite{CS, CS2}. 
In particular, $\X''_{\rm min,0}$ 
(together with the corresponding negative valuation) 
is K-polystable in the sense of \cite[\S 5]{Song}.

\item \label{Kss.Kps}
There is an affine ($\Q$-)Gorenstein test configuration of $\X_{\rm min,0}'$ degenerating to $\X_{\rm min,0}''$. 

\end{enumerate}
\end{Thm}

We prove Theorem \ref{bblcst:def} 
in the subsection \S \ref{sec:AGproof}. 

\begin{Rem}[Generalization to klt pairs]\label{ll.rem}
Generalizations of the above results to general non-Gorenstein log terminal morphisms, and even to 
log pair setup, should be straightforward adaptation 
which we omit here to avoid expository complication. 
It would involve  
a subtle construction of $\Q$-Gorenstein
\footnote{``locally stable" in the sense of \cite{Kol23}}
parameter space and its routine logarithmic extension, 
because of the fiberwise non-Gorensteinness. 
It should be solvable by the use of universal hull construction 
\cite{Kol08}, \cite[Chapter 9, \S 9.9]{Kol23}, 
after the flattening stratification. 
\end{Rem}

For the sake of organization of our explanations, 
here we focus and finish the issue of finiteness of bubbles, 
for fixed $\pi$, 
as a part of the previous theorem \ref{bblcst:def}. 
We use \cite{XZ24} but we note that 
a differential geometric counterpart 
is proven before for the case when $\X_t (t\neq 0)$ are all smooth, 
compactifiable to  
metrized K\"ahler manifolds, 
in \cite[Corollary 1.12]{Song}. 

\begin{Thm}\label{ssr.sing}
We fix a positive integer $n$. 
In the setup \ref{Setup}, \ref{Setup2}, 
we modify $\X$ by $\X''_{\rm min}$ as in Theorem \ref{bblcst:def}. 
Then this modification process does not continue infinitely many times for 
$\dim(X)=n$. Hence, we algebraically obtain a 
sequence of (weighted) algebraic bubbles 
$Z_0,\cdots,Z_{2k}$ similarly to \cite[Cor 1.12]{Song}. 

In particular, as the terminal object of this repetition, 
possibly after a finite base change of $\Delta$, 
there is an affine faithfully flat klt Gorenstein 
family $\mathcal{Y}\to \Delta$ 
which coincides with $\mathcal{X}\to \Delta$ over $\Delta\setminus \{0\}$, 
whose central fiber is $\mathcal{Y}_0=Z_{2k}$, and such that the section $\sigma$ 
extends to a section $\tau\colon\Delta\to\mathcal{Y}$ satisfying 
\begin{align}
\label{ssr1}\widehat{\vol}_{\Y_0}(\tau(0))&=\widehat{\vol}_{\Y_t}(\tau(t))\\ 
\label{ssr2}&=\widehat{\vol}_{\X_t}(\sigma(t))
\end{align}
for every $0<|t|\ll 1$. 

\end{Thm}

\begin{proof}
At every nontrivial repetition of Theorem \ref{bblcst:def}, the normalized volume 
at the marked point of the central fiber strictly increases.  It remains bounded above by 
$n^n$, while \cite[Theorem 1.2]{XZ24} shows that, in fixed dimension, the set of local 
volumes has no positive accumulation point.  Hence the process terminates after finitely many steps. 
We take $\mathcal{Y}$ to be the terminal family.  The first equality in \eqref{ssr1} holds 
because otherwise Theorem \ref{bblcst:def} could be applied once more, and \eqref{ssr2} follows 
from the identification of $\X$ and $\Y$ over $\Delta\setminus\{0\}$.
\end{proof}

Note that the last statement of the above Theorem 
\ref{ssr.sing} 
is a semistable reduction type theorem corresponding to virtual compactness of 
 ``moduli" of singularities with fixed normalized volumes. 
In particular, 
the following corollary to above, of a weak version of the 
resolution of singularities type, 
is pointed out by M.Mauri. 

\begin{cor}[As a variant of alteration] 
Consider any affine pointed flat morphism 
$\pi\colon \X\to \Delta\ni 0$ of relative dimension $n$, 
with a section $\sigma\colon 
\Delta\to \X$ such that fibers $\pi^{-1}(s)=\X_s\ni \sigma(s)$ 
are all smooth points (resp., log terminal) 
for $s\neq 0$ (resp., $s=0$), as in Setup \ref{Setup}. 

Then there is a birational modification 
$\tilde{\pi}\colon 
\Y\twoheadrightarrow \Delta$ along $\X_0$, possibly after a 
finite base change of $\Delta$, such that 
$\sigma|_{\Delta\setminus 0}$ 
extends to whole $\Delta$ as $\sigma_{\Y}\colon \Delta\to \Y$, 
$\tilde{\pi}$ 
is smooth surjective around $\sigma_{\Y}(s)$ for any $s$ 
near $0$. 

If $n=2$,  or if $\pi$ is a log terminal morphism more generally, 
one can construct $\Y$ with $\sigma_{\Y}$ which are 
\'etale locally 
trivial around $\sigma_{\Y}(0)$. 
\end{cor}

Our main reason for writing this is to present  our different method 
than usual resolution of singularities, 
rather than the result itself. 

\begin{proof}
Suppose the last assertion is not satisfied for 
the original $\X$. Otherwise, the statement is obvious. 
Now we apply 
Theorem \ref{bblcst:def} and Theorem \ref{ssr.sing} 
to obtain the deepest bubbling family 
$\tilde{\pi}$. 
From \eqref{ssr1} and \eqref{ssr2} of 
Theorem \ref{ssr.sing}, $\widehat{\rm vol}_{\Y_0}
(\sigma_{\Y}(0))=n^n$. Then, \cite[Theorem A.4]{LiuXu}
 (see also implicit arguments in \cite{HS, SS, LX}) 
 shows that $\sigma_{\Y}(0)\in \Y_0$ is a smooth point. 
 Hence we complete the proof of the former claim for general $n$. 
 The claim for $n=2$ follows from the same arguments, 
 once we combine with the simpler proof of Theorem 3.6 of 
 \cite{Od24a} using the comparison of the orders of the local 
 fundamental groups. 
\end{proof}

Further, note that 
the above analog to the resolution of singularities 
$\tilde{X}$ is often 
taken canonically if $\X_0$ has a quasi-regular 
minimizer of $\widehat{\rm vol}(-)$. 
 
In our companion paper \cite{Od24a}, 
we restrict our focus to cones where the similar method applies to 
prove the ``proper"-ness part of the following theorem. 

\begin{Thm}
[K-moduli of Fano cone \cite{Od24a}]\label{Kmod.sing}
Consider the set (resp., all $\Q$-Gorenstein families) of 
the isomorphic classes of 
K-polystable (resp., K-semistable) $\Q$-Fano cone of 
$n$-dimension, with 
the fixed normalized volume $V$ of the vertices. 
Then, 
it forms a good proper (and separated) 
moduli algebraic space over $\mathbb{C}$ 
and its stacky enhancement, 
generalizing the K-moduli of $\Q$-Fano varieties. 
\end{Thm}

\subsubsection{Comparison with 
differential geometric bubbles}

As a link to differential geometric bubbling  
i.e., rescaled up pointed Gromov-Hausdorff limits, 
which are studied in 
\cite{Song, dBS}, 
we have the following natural question. 

\begin{Ques}[Identification with differential geometric bubbling]
Does our $\X''_{\rm min,0}$ 
in Theorem \ref{bblcst:def}, 
applied with some $\mathcal{T}''$ 
and good ambient local coordinates, 
coincide with the (differential geometric) 
minimal bubble in the 
sense of \cite[Definition 1.9]{Song} (Definition \ref{DG.min}), 
obtained in \cite[Theorem 1.10]{Song}? 
\end{Ques} 

Note that we refer to a priori 
non-canonicity of the constructions in Theorem \ref{bblcst:def}, 
which we expect to parallel the choice of the family of 
K\"ahler metrics $\{g_t\}_t$ on 
$\X$. In the following $1$- or $2$-dimensional classical cases 
considered in \cite{dBS}, the algebraic and differential-geometric 
constructions can be compared directly. 
Below, we denote an analytic neighborhood of $0\in \Delta$ 
as $\Delta_t$. 

\begin{Thm}[Confirmation - identification with differential geometric 
bubbling]\label{AGDG.cl}
The algebraic construction of Theorem \ref{bblcst:def} coincides 
with the differential-geometric bubbling construction in at least the 
following three cases. These form examples of Setup \ref{Setup}. More 
precise statements are itemized as follows.

\begin{enumerate}
\item \label{cl1} (``Log curve" case)
\footnote{To be precise, our Theorem \ref{bblcst:def} does not directly 
apply to this log setup but the proof applies easily in this case (cf., Remark \ref{ll.rem}). We clarify the meaning during the proof.} 
If we consider $\mathcal{X}:=\C\times \Delta_t\to \Delta_t$  
with flat (singular) K\"ahler metrics 
with conical singularities $p_i(t)\in \C$ for 
$i=0,\cdots,k$ of angle $2\pi(1-\beta_i)$ with $0<\beta_i<1$, 
the sequence of bubblings constructed algebro-geometrically 
by Theorem \ref{bblcst:def} 
(cf., Remark \ref{ll.rem}) 
can be taken canonically, 
coincides with 
that of \cite[\S 2.2, Lemma2]{dBS} and hence 
coincides with differential geometric bubbles i.e., 
rescaled pointed Gromov-Hausdorff limits 
(or, polarized limit space in the sense of \cite{DSII}, 
to be more precise). 

\item \label{cl2} ($A_k$-type ALE case ($n=2$))
Consider general smoothing affine 
family of ADE singularity 
$\mathcal{X}(\supset \Delta_t)\to \Delta_t$ where 
$X=\X_0$ has only one singularity at  $\sigma(0)=0$ 
which is of $A_k$ type. 
Suppose we have 
a continuous family of 
ALE hyperK\"ahler (orbi-)metrics $g_t$ on $\X_t$ 
whose K\"ahler forms 
$(\omega_{I,t},\omega_{J,t},\omega_{K,t})$ 
($I$ is the original complex structure for $X_t$s 
in Setup\ref{Setup}) 
satisfy that 
\begin{itemize}
\item  
$\omega_{I,t}$ is $\partial \overline{\partial}$-exact 
\item 
$\omega_{J,t}+\sqrt{-1}\omega_{K,t}$ is holomorphic 
on $\X\setminus X$. 
\end{itemize} 
Then, 
the sequence of bubblings constructed algebro-geometrically 
by Theorem \ref{bblcst:def} 
can be taken canonically again, 
coincides with 
the construction of \cite[\S 3.3, Theorem 3]{dBS}, and hence 
matches to differential geometric bubbles. 

\item \label{cl3} ($A_1$ case with any $n$)
If $\mathcal{X}\to \Delta_t$ 
is open subset of proper holomorphic family over $\Delta$ 
which is the smoothing of 
higher dimensional $A_1$-singularity, 
the minimal bubbling constructed algebro-geometrically 
by Theorem \ref{bblcst:def} is 
$V(z_0^2+\cdots+z_{n}^2-1)\subset \A^{n+1}$ 
and coincides 
with the differential geometric minimal bubble with the 
Stenzel asymptotic conical metric (\cite{EH, Ste}). 
Here, $V(-)$ denotes the zero locus (variety) defined by the equation. 

\end{enumerate}
\end{Thm}

Also, as a related result to the above \eqref{cl2} and \eqref{cl3}, 
after our paper, 
\cite{Tazoe} showed similar statements for the degeneration of 
the hyperK\"ahler metrics on type I degenerations of polarized K3 surfaces. 

Henceforth, 
we proceed to the proofs of the above 
Theorems \ref{bblcst:def} and \ref{AGDG.cl}. 

\subsubsection{Details of the algebraic 
algorithmic construction}
\label{sec:AGproof}

In this subsubsection, we prove 
Theorem \ref{bblcst:def}. 

\begin{proof}[proof of Theorem \ref{bblcst:def}]
First we discuss the construction of 
$\X'_{\rm min}$, its central fiber $\X'_{{\rm min},0}$ 
with its degeneration to a K-semistable klt cone $W\curvearrowleft T$ (recall at \S \ref{sec:rev}). 
We consider the deformations of $W$ and here we a priori face the problem of infinite-dimensionality of the deformation space of affine varieties, 
especially those with non-isolated singularities 
(see \cite[\S 6, 6.9]{SZ}). 
Nevertheless, following the idea of negative weight deformation (\cite{CH}, 
also cf., \cite{Song, dBS}), 
we can circumvent the problem. 

First, we take a temporary embedding of $W$ to $\A^l$ 
via (non-constant) $T$-semiinvariants $z_1,\cdots,z_l$ so that the vertex 
maps to the origin. Suppose its Reeb vector field  $\xi$ 
which corresponds to the volume minimizer and has 
weights $w_1,\cdots,w_l \in \R_{>0}$ for each coordinate of $\A_{\mathbb{C}}^l$ i.e., 
${\rm deg}_{\xi}(z_i)=w_i$ (cf. \cite{DSII, 
Od24a, Od24b}). 
We consider finite generators $f_1,\cdots,f_N$ of the 
defining ideal which gives 
$W\subset \A_{\mathbb{C}}^l$ and set $d_i:={\rm deg}_{\xi}(f_i)$. 
Consider all affine deformations of $W$ 
of the $\xi$-negative type 
$V(\{f_i+h_i\}_{i=1,\cdots,N}\mid 
0<{\rm deg}_{\xi}(h_i)<{\rm deg}_{\xi}(f_i))$. 
Also compare more intrinsic but conditional related 
construction: Corollary \ref{Kur}. 
Then, by applying the flattening stratification 
we obtain a $T$-equivariant 
affine flat 
family $\pi_{\mathcal{Y}}\colon 
\mathcal{Y}\to {\rm Def}^{-}(W)\subset \A^m$ for $m\gg 0$. 
By \cite[Proposition 4.1]{NO}, $W$ is Gorenstein. Further, 
$\pi_\Y$ is a Gorenstein morphism in the relative sense  
(and hence locally stable in the sense of \cite{Kol23}). 
Indeed, the relative Gorenstein locus in ${\rm Def}^-(W)$ is open, 
contains the origin, and invariant under 
the definite $T$-action, hence it is automatically the whole 
${\rm Def}^-(W)$. 
We consider its $T$-invariant 
locally closed subset $$Z^+:=\{b\in {\rm Def}^-(W)\mid 
\widehat{\rm vol}(0\in \pi_{\Y}^{-1}(b))=\widehat{\rm vol}(0\in W))\},$$ 
by \cite{BlumLiu}. 
By replacing ${\rm Def}^{-}(W)$ by 
a $T$-invariant open affine neighborhood by Sumihiro's theorem (\cite{Sumihiro}), 
we can suppose $Z^+$ is its closed subset. 
Further, by \cite{Chen}, $\mathcal{O}_{\pi_{\Y}^{-1}(Z^+)}$ 
has ideal sequences $\{I_{\lambda_i}\}_{\lambda_i\in S(W)\subset \R_{>0}}$ 
such that for any $b\in Z^+$ 
$\{I_{\lambda_i}|_{\pi_{\Y}^{-1}(b)}\}_{\lambda_i\in S(W)\subset \R_{>0}}$ 
gives the degeneration to the K-semistable Fano cone 
\`{a} la \cite{DSII}. 
Here, $S(W)$ is the holomorphic spectrum of $W$ also 
in the sense of \cite{DSII}. 
Recall these at \S \ref{sec:rev}. 

Now, we take a positive integer $m$ and 
shrink ${\rm Def}^-(W)(\ni [W])$ enough again 
to a $T$-invariant open affine neighborhood 
such that 
all $i=1,\cdots,m$, 
$I_{\lambda_{i}}/I_{\lambda_{i+1}}$ are all free $\O_{Z^+}$-modules of 
rank $r_i$, 
and have $T$-equivariant splittings of 
$I_{\lambda_i}\twoheadrightarrow I_{\lambda_{i}}/I_{\lambda_{i+1}}\simeq 
\O_{Z^+}^{r_i}$ 
which corresponds to $T$-semiinvariant 
$\{z_{i,j}\in I_{\lambda_i}\mid 1\le i\le m, 
1\le j\le r_i\}$ 
(\cite[2.5]{IltenSuss.vb}). 
By the complete reductivity of $T$, one can also extend 
$z_{i,j}$ to $T$-semiinvariant regular functions on whole $\Y$, 
which vanishes at the $0$-section. 
With large enough $m$ and shrunk enough ($T$-invariant) 
$[W]\in {\rm Def}^-(W)$, 
we can and do assume the above 
$z_{i,j}$ generates $\Gamma(\mathcal{O}_{\pi_{\mathcal{Y}}^{-1}(b)})$ 
as $\mathbb{C}$-algebras 
for any $b\in  {\rm Def}^-(W)$. 
From the construction, it follows that 
for any $b\in Z^+$, the natural $T$-action with the weights 
$w_i$ give the degeneration to K-semistable Fano cone. 
We here reset $l:=\sum_{i=1}^m r_{i}$ and $\xi=(w_1,\cdots,w_1,\cdots,
w_m,\cdots,w_m)$. Then, the above $z_{i,j}$ 
gives a large reembedding of $\mathcal{Y}\hookrightarrow \A_{\mathbb{C}}^l\times {\rm Def}^-(W)$.  

Next, we enlarge ${\rm Def}^-(W)$ a little just to include $\X_t$s. 
Suppose the reembedding $W\subset \A_{\mathbb{C}}^l$ constructed above is defined by $c$ 
polynomial equations $g_1,\cdots,g_c$. 
Then, take an arbitrary extended embedding of $\X$ to 
$\A_{\mathbb{C}}^l\times_{\mathbb{C}} \Delta$ and suppose that 
$\X_t$ for general $t$ 
is defined by equations of $\xi$-degrees bounded above by some $C>0$. 
Now we consider the set of all possible polynomials $g_1,\cdots,g_c$ on $\A_{\mathbb{C}}^l$ of the degrees at most by $C$, without constant terms, 
and consider the corresponding universal family 
$\mathcal{U}\twoheadrightarrow 
\A_{\mathbb{C}}^p$ 
for $p\gg 0$. There is a corresponding morphism 
$\iota\colon \Delta\to \A_{\mathbb{C}}^p$ to this $\X\to\Delta$. 


Then we obtain $\overline{T\cdot \iota(\Delta)}(\subset \A_{\mathbb{C}}^p)$, apply the flattening stratification, 
and take the strata $H$ which includes 
$[W]$ and intersect with $\iota(\Delta)$. 
From the construction, $H$ is automatically $T$-invariant and 
the obtained affine family over $H$ is flat, 
which we denote by $\mathcal{V} \to H$. 
(Recall that its restriction to ${\rm Def}^-(W)$ is a Gorenstein family, 
by the preceding argument.) 
This $H$ will be our parameter space to work on. 

As $\Delta$ in Setup \ref{Setup} 
is a pointed smooth affine $\mathbb{C}$-curve 
so that $H$ underlies a structure of algebraic 
$\mathbb{C}$-scheme. 
Here, 
this $T$-action on 
$\iota(\Delta)\subset \A_{\mathbb{C}}^p$ is given via 
the re-embedding coordinates and is the same for the 
whole $T$-action on $M$. 

We denote 
$\mathcal{M}:=[H/T]$. 
(For those who prefer rank $1$ setup, one could alternatively take 
$\mathcal{T}$ a (non-canonical) affine test configurations 
induced by 
an algebraic subtorus $T'$ of $T$ with the rank $1$ 
(the symbol $\mathcal{T}'$ is reserved for different 
family to appear later as it is more natural notation for it), 
as a subfamily of the whole family $\X_\tau\to U_\tau$ 
of \cite[Example 2.9]{Od24b}. 
For $r>1$ case, this $T'$ is clearly 
non-unique as a rational approximation of the 
algebraic positive degeneration $X\rightsquigarrow W$ 
(\cite[\S 3]{LX}, \cite[\S 2]{Od24b}).) 

From our construction, the resulting $\varphi\colon [\mathcal{V}/T] \to \M$ is a 
($\Q$-)Gorenstein log terminal family over an Artin stack and 
all the fibers contain the origin $0$. 

Now it follows that $[Z^+/T]$ 
has a structure of 
canonical higher $\Theta$-strata 
(\cite{Od24b}) 
on $\mathcal{M}$ near $[W]$ 
by the above construction. 
This follows from the above construction of the 
$T$-acted base, using \cite{Chen}, 
together with the following fact: 
if you have a log terminal affine morphism with base 
points 
$f\colon \tilde{\X}(\supset \sigma(S)=\{0\}\times S)\to S$ and 
a uniform positive weight degeneration 
$\tilde{X}_s\ni 
\vec{0} \rightsquigarrow W_s \ni \vec{0}$ 
of the pointed klt affine variety $X_s\ni 
\vec{0}$ to 
its K-semistable cone $W_s$, 
i.e., if 
there is a flat family 
over $S\times_\mathbb{C} U_\tau$ with a connected algebraic $\mathbb{C}$-scheme $S$ 
and an affine toric variety $U_\tau$ which corresponds to 
a rational polyhedral cone $\tau$, 
the corresponding minimal 
normalized volume $\widehat{\rm vol}(f^{-1}(s)\ni 0)$ 
is constant i.e., not depending on $s\in S$. 
This follows from 
the comparison with the leading coefficient  of the index character  
(cf., \cite{MSY, Li, CS} and also \cite[4.1]{BB}, \cite[\S 2]{Od24a}). On the other hand, 
over a neighborhood of $[W]\in Z^+$, 
our new $T$-action gives a 
higher $\Theta$-strata by \cite{Chen}. 

Hence, we can apply semistable reduction type 
theorem \cite[Theorem 3.5 (also cf., 1.1)]{Od24b} 
(which generalizes \cite{Langton, AHLH}) which gives 
a Langton type 
modification of the original morphism 
$f\colon \Delta \to [H/T]$ (or its formal germ 
${\rm Spec}\C[[t]]\to [H/T]$) that corresponds to $\X$. 
As its consequence, for some 
finite covering $\Delta^{(d)}\to \Delta$  for some $d\in \Z_{>0}$ 
(see Setup \ref{Setup}), 
we obtain $\Delta^{(d)}\to H$ which gives an affine flat family. 
We are now going to 
confirm the required conditions of $\X'_{\rm min}$
in our Theorem \ref{bblcst:def}, and confirm that 
$\mathcal{W}$ of \eqref{Kss.Song} 
also naturally show up during the construction. 


Since $\X'_{\rm min,0}$ is a negative weight deformation of $W$, hence Gorenstein 
by the preceding argument, $\X'_{\rm min}\to \Delta^{(d)}$ is a Gorenstein morphism. 
The K-semistability of $\X'_{\rm min,0}$ in the sense of \cite[\S 5]{Song} 
follows from the construction as $W$ is a K-semistable log terminal cone 
($\Q$-Fano cone) 
in the sense of \cite{CS, CS2}. Furthermore, from the positivity of the 
weights on the coordinates of $T'$-action and the definition of elementary modification 
\cite[\S 3]{AHLH}, we conclude that the obtained test configuration 
$\X'_{\rm min,0} \rightsquigarrow W$ is a negative weight deformation in the sense of 
\cite{CH} and \cite[\S 2]{Od24b} of rank $1$. 
A rank-one subtorus $T'$ is generally non-unique when the original direction is irrational. 
The independence of the resulting central fiber follows instead from the higher-rank 
canonicity statement \cite[Theorem 3.5(ii)]{Od24b}.  This is why we use 
the higher-rank theorem \cite[Theorem 3.5]{Od24b}, rather than only the original 
rank-one result \cite[\S 6]{AHLH}.


\begin{figure}[ht]
\centering
\(\vcenter{\hbox{\xy 0;<0.9mm,0mm>:
  (0,0)*{};(100,76)*{};
(78,66)*+{\ensuremath{\mathcal X}};
(78,61);(78,52) **\dir{-} ?>*\dir{>};
(78,47)*+{\ensuremath{\Delta}};
  (20,45)*+{\ensuremath{\mathcal X_{\min}^{\prime}}};
  (20,40);(20,33) **\dir{-} ?>*\dir{>};
  (20,28)*+{\ensuremath{\Delta^{(d)}}};
(18,10);(55,10) **\dir{-} ?>*\dir{>};
(89,19);(59,11) **\dir{--} ?>*\dir{>};
  (57,10)*{\ensuremath{\bullet}};
  (22,5)*+{\ensuremath{\X'_{\rm min,0}}};
  (57,5)*+{\ensuremath{W}};
  (89,19)*{\ensuremath{\bullet}};
  (93,15)*+{\ensuremath{X}};
  (40,14)*+{\ensuremath{\mathcal W}};
(56,69);(89,19) **\crv{~**\dir{--}(69,64)&(75,38)};
  (49,68);(18,10)
    **\crv{~**\dir{--}(48,52)&(31,25)}
    ?>*\dir{>};
\endxy}}\)
\caption{The construction of $\X'_{\rm min,0}$}
\label{fig:semistable-minimal-bubble}
\end{figure}


Next, we discuss 
algebraic construction of 
$\X''_{\rm min}$ from $\X'_{\rm min}$, 
which follows a similar line of ideas 
as those for 
$\X'_{\rm min}$, 
but also with some differences. 
We take a $T$-equivariant test configuration of 
$W$ to the K-polystable degeneration $C$ as in \cite[Theorem 1.2]{LWX}. 
Recall that $C$ is unique while the total space of the family is not. 
Denote the acting $\G_m$ as $T''$ and corresponding base i.e., the 
partial compactification of $T''$ as 
$\overline{T}''$ which is just $\A^1$. This notation is introduced to avoid 
confusion with other algebraic tori appeared above. 
Recall that the construction in {\it loc.cit} is far from canonical. 
We still take one of them and denote the corresponding 
test configuration as $\mathcal{T}''\to \overline{T}''$. 
If we write $\mathcal{T}''={\rm Spec}_{\C[t]}(\oplus_{e\in \Z, \vec{m}\in M}
R^{(\vec{m})}_e)$ for an increasing filtration $R^{(\vec{m})}_e\subset 
R^{(\vec{m})}_e$ on the $\vec{m}$-eigenspace of $\Gamma(\mathcal{O}_W)$, 
and take finite generators of the 
double graded ring $(\oplus_{e\in \Z, \vec{m}\in M}
R^{(\vec{m})}_e)$. It allows to reembed $\mathcal{T}''\hookrightarrow \A^l\times 
\A^1$ in a $T\times \G_m$-equivariant manner and we similarly consider the 
generators of $\Gamma(\mathcal{O}_{\X'_{\min}})$ as $\Gamma(\mathcal{O}_\Delta)$-
algebra. 
After adjoining the test configuration $\mathcal{T}''; W\rightsquigarrow C$ 
and the scale down degeneration of $\X'_{\rm min,0}$ to $W$, to 
a new bounded-degree coefficients parameter affine space as before, we apply 
the ordinary flattening stratification once more. 
As a consequence, 
we can and do assume $\mathcal{T}''$, $\X'_{\rm min}\rightsquigarrow W$, 
$\mathcal{X}'_{\rm min}$ are all realized over the same base $H'$ 
as ($\Q$-)Gorenstein family as before. 
Now we consider a couple of morphisms 
$\Spec \C[[t]] \to H$ 
corresponding to 
$\mathcal{T}'$ of degeneration of $X' \rightsquigarrow W$, and corresponding to $\mathcal{T}''$ 
whose degeneration is $W \rightsquigarrow C$. 

We denote the base, with its associated moduli map, of the 
negative weight degeneration $\X'_{\rm min,0}\rightsquigarrow W$ 
by $\Delta'\to H'$. Then, consider the natural rational map from 
$m\colon \Delta'\times \overline{T}'' \to H'$. 
If it extends to a morphism, one can take the image of 
$c\times \{0\}(\in \Delta'\times \overline{T}'')$ for $c\neq 0$ 
as $\X''_{\rm min,0}$. Otherwise, we apply 
\cite[Theorem 6.3]{AHLH} (or \cite[Theorem 3.5]{Od24b}) to them, 
with respect to the simple 
$\Theta$-strata $[\overline{m(\Delta'\times T'')}/T'']$ supported at 
$m(0\in \Delta)$. 
Here, $\overline{m(\Delta'\times T'')}$ means the closure of 
the image $\overline{m(\Delta'\times T'')}$ of $m$. 
Then, it follows that 
possibly after passing to a finite ramified covering 
$\Delta^{(d)}\to \Delta$, 
we obtain an isotrivial family which degenerates $\X'_{\rm min,0}$ 
to $\X''_{\rm min,0}\subset \A_{\mathbb{C}}^l$. 
In this case, we can avoid the full use of \cite{AHLH, Od} 
by a more standard direct argument using a toric surface, as follows. 

Note that $T'$ and $T''$ are commutative so that they are both algebraic subgroups of 
an algebraic torus $T_2$ of rank $2$, 
which acts equivariantly on the family 
$\mathcal{U}\to \A_{\mathbb{C}}^m$ and its base change  
$\mathcal{V}\to H''$. 
Since the obtained elementary modification maps 
$(0,0)\in [\A^2/\G_m]={\rm ST}_{k[[t]]}$ (cf., \cite{HL, AHLH}) 
to 
$[C\in \mathcal{M}]$ with the inertia group $T$, 
the elementary modification is represented by a 
$T_2$-equivariant morphism $f\colon \A_{\mathbb{C}}^2\to H''$. The normalization of its 
image in $H''$ is a toric surface $T_{N\simeq \Z^2}\Sigma$ where 
the fan contains a ray $\R_{>0}(1,0)$ (resp., $\R_{>0}(0,1)$) which corresponds to 
$\mathcal{T}'$ (resp., $\mathcal{T}''$) and its 
support ${\rm Supp}(\Sigma)=\R_{\ge 0}^2$ from the construction. 
Take the ray $\R_{>0}(s_1,s_2) (s_i\in \Z_{>0})$ which is adjacent to $\R_{>0}(1,0)$ 
and consider the one parameter subgroup $\lambda \colon \G_m\to T_2$. 
Then it gives a test configuration $\X'_{\rm min,0} \rightsquigarrow \X''_{\rm min,0}$ for some $X''\in S$ such that 
the inverse action of $\mathcal{T}''$ degenerates $\X''_{\rm min,0}$ to 
$C$ 
as desired.

Now we glue a family $\X'_{\rm min,0}$ and $\pi_{X'}$, by the same technique as 
\cite[\S 3]{LWX}, \cite[\S 3]{ops}, 
to obtain the desired $\X''_{\rm min,0}\to \Delta$, as a vertical birational modification of $\X'_{\rm min,0}\to \Delta$, 
which satisfies the conditions listed in the statements of our theorem. 
The only part which does not follow immediately, is the equality assertion of \eqref{vol.ineq}. It holds because of 
the condition satisfied by the ``minimally" constructed 
$\X'_{\rm min,0}\ni \sigma'(0)$ i.e., that it achieves the possible minimum 
normalized volume, combined with the semicontinuity \cite[Theorem 1]{BlumLiu}, 
\cite[Theorem 1.3]{Xu} of $\widehat{\rm vol}$ 
applied to $\pi_{X'}$. 
Here we complete the construction of $\X''_{\rm min,0}$ from 
$\X'_{\rm min,0}$. If $n=2$, since any log terminal singularities 
are quotient singularities so that K-semistable Fano cones are 
automatically K-polystable, we have $\X''_{\rm min,0}=\X'_{\rm min,0}$ from 
our construction. 
We complete the proof. 
\end{proof}


\begin{Rem}[Relation with deformation theory]
Note that deformation theory of singularity usually deals with 
deformations of {\it germs} so that a priori it can not be applied directly 
to the construction of these bubbles - a certainly global object (affine 
variety, not just a germ). 

Nevertheless, 
under certain situation, we can compare with deformation theory of $W$ 
which gives more canonical expression. 
For instance, suppose 
$W$ has only isolated singularity 
and has unobstructed deformation i.e., 
with smooth 
$T$-equivariant deformation space 
$[W]\in {\rm Def}(W)$ which is algebraized by 
\cite{Artin, Pink}. 
Then, by \cite{KR, BH, LP}, 
$[W]\in {\rm Def}(W)$ is 
$T$-equivariantly 
isomorphic to 
$0\in \A_{\mathbb{C}}^d$ with linear $T$-action.  
This gives a $T$-equivariant decomposition 
${\rm Def}(W)=\oplus_{i\in I} 
V_i$ to $T$-eigenspaces $V_i$. 
We denote the $\xi$-weights on $V_i$ as $W_i$. 
Now we add an assumption that 
$\iota \colon \Delta\to {\rm Def}(W)$ which corresponds to 
$\X$ has non-vanishing differential at each 
$V_i$: we call this condition {\it 
weakly non-degenerating}, comparing with 
stronger conditions in \cite{BR}. A subtle 
point is that because of a priori non-universality of 
${\rm Def}(W)$, $\iota$ is not canonical 
and only its differential $\partial \iota_{t=0}$ 
is canonical by the versality. 
On the other hand, 
the construction of $\X'_{\rm min}$ 
shows that they are determined by 
$\partial_t \iota|_{V_i}$ 
for those $V_i$ with $W_i=1$. 
Thus, the weak non-degeneracy assumption implies 
that the above construction of 
$\X'_{\rm min,0}$ is canonical. See also Tazoe's work \cite{Tazoe} on the case of K3 surfaces degenerating with ADE singularities.  
\end{Rem}

\subsubsection{Comparison with the differential geometric construction}

Now we prove Theorem \ref{AGDG.cl}. 
The main point is confirmation of 
the compatibility of our construction 
in the proof of Theorem \ref{bblcst:def}  with 
the case by case analysis in 
\cite{dBS}, in each case. 

\begin{proof}[proof of Theorem \ref{AGDG.cl}]
In all the three cases of  Theorem \ref{AGDG.cl}, 
$X\ni \sigma(0)$ is locally an 
isolated singularity with local $\C^*$-action 
(with conical singularities) 
and hence in particular we can consider 
its finite dimensional smooth (unobstructed) 
deformation space ${\rm Def}(X\ni \sigma(0))$ 
of the singularity. 

First, we consider Case \eqref{cl1}. 
In this case, 
we consider deformations of 
$X\ni \sigma(0)$ as 
families of 
(at most) $k+1$-pointed $\P^1$. Its 
semi-universal deformation is 
straightforward (cf., also \cite[\S 2.2, 
\S 2.3]{dBS}). 
Indeed, 
there is a finite covering 
(the Weyl type covering) 
$\D:=\C^{k+1}\to {\rm Def}(X\ni \sigma(0))=\C^{k+1}/W\simeq 
\C^{k+1}$,  
where $W\subset \mathfrak{S}_{k+1}$ is the 
subgroup of the 
symmetric group of degree $k+1$, 
as the symmetry of the tuple 
$(\beta_0,\cdots,\beta_k)$, 
which acts by the permutation. The quotient is a Galois  covering with the Galois group $G$. 
We consider the morphism $\Delta_t\to {\rm Def}(X\ni \sigma(0))$ 
which is associated to $\X(\supset \Delta_t)\twoheadrightarrow \Delta_t$ 
and its lift $\Delta_t\to \D$ possibly after 
replacement of $\Delta_t$ by its 
finite covering (a component of the base change 
$\D\times_{{\rm Def}(X\ni \sigma(0))}\Delta_t$). 
We still denote the finite covering as $\Delta_t$ for simplicity, 
and denote as $\varphi\colon \Delta_t\to \D=\C^k$, 
$t\mapsto (s_0(t),\cdots,s_k(t))$ 
where each $s_i(t)$ is holomorphic with $s_i(0)=0$. 
In particular, the vanishing order of $s_i(t)$ at $t=0$ 
is positive which we denote as $d_i$. We set 
$d:=\min_i \{d_i\}$. 

Note that the induced action of $\C^*$ as a subgroup of the 
automorphism group of $X\ni \sigma(0)$ on 
$\mathcal{D}$ is with the weights $(-1,\cdots,-1)$. 
Thus, from the proof of Theorem \eqref{bblcst:def}, 
$\X'_{\rm min}$ naturally 
corresponds to 
\begin{align}\label{8}
t\mapsto \varphi'_{\rm min}(t):=(t^{-d}s_0(t),\cdots,t^{-d}s_k(t))\in \C^{k+1}=\D.
\end{align}
Indeed, for $\epsilon>0$, 
we have that $\lim_{t\to 0} 
t^\epsilon \varphi'_{\rm min}(t)
=[X\ni \sigma(0)=0]$, 
while $\varphi'_{\rm min}(0)\neq 0$. 
On the other hand, it is easy to see that the flat (conical 
singular) K\"ahler metric on $\C$ is uniquely 
described as 
$$\prod_{0\le i\le k}|z-s_i(t)|^{\beta_i-1}|dz|$$ 
(compare \cite[Lemma 1]{dBS}), it follows that the 
differential geometric bubbling is obtained as 
the flat singular $\C$ given by $\varphi'_{\rm min}(0)$. 
By replacing $\varphi$ by $\varphi'_{\rm min}$ 
and iterating the procedure, we complete the proof of 
\eqref{cl1} of Theorem \ref{AGDG.cl}. 

The proof for the case \eqref{cl2} is 
more or less parallel to the case \eqref{cl1} as above. 
We consider the well-known 
miniversal deformation space of $A_k$-singularity 
by Brieskorn-Grothendieck-Kas-Schlessinger-Slodowy. 
Recall it is also reconstructed as 
Kronheimer's construction of 
a family of 
hyperK\"ahler quotients \cite{KronI}, 
which we use. The deformation space can be taken as $\C^k$ 
and there is a natural Weyl cover $\mathcal{D}$ 
which kills the monodromy. 

In our setup, due to the existence of ALE 
hyperK\"ahler metrics $g_t$ with the 
$\partial \overline{\partial}$-exact K\"ahler form 
$\omega_{I,t}$, $\X\to \Delta_t$ 
can be written as (after base change to the Weyl cover i.e., 
take a component of $\Delta_t\times_{{\rm Def}(X\ni \sigma(0))}\D$) 
$$X_t=V(xy-\prod_{i=0}^k (z-s_i(t))\subset \A_{\C}^3,$$ 
where $s_i(t)$ is holomorphic. 
Indeed, each $X_t$ can be compared with 
the Gibbons-Hawking ansatz for appropriate 
$(0,s_i(t))\in \R\times \C \simeq \R^3 (i=0,\cdots,k)$ 
with the harmonic potential 
\begin{align*}
   \sum_{0\le i\le k} \frac{1}{|z-s_i(t)|} 
\end{align*}
and the corresponding connection form as an essentially unique 
solution to 
the Bogomolny-type monopole equation, 
as explained in \cite[p.302-303]{LeBrun}, \cite[\S 3.2, \S 3.3]{dBS} and it follows 
from Kronheimer's Torelli-type 
theorem \cite{KronII} 
that they are isomorphic as hyperK\"ahler varieties. 
It is obvious that $s_i(t)$ should be continuous but 
is even holomorphic 
by the second bulleted assumption on the holomorphicity of 
$\omega_{J,t}+\sqrt{-1}\omega_{K,t}$ on $\X\setminus X$. 
Now we can apply the discussion very similar to the proof of 
case \eqref{cl1} i.e., consider 
the vanishing order of $s_i(t)$ at $t=0$ which 
is positive and denote by $d_i$ and set 
$d:=\min_i \{d_i\}$. Then we consider the map 
to 
the Weyl cover 
$\varphi'_{\rm min}\colon 
\D(=\C^{k+1})\to {\rm Def}(A_k)$ exactly as 
\eqref{8}. Since the induced action of $\C^*$ as a subgroup of the 
automorphism group of $A_k$-singularity 
$X\ni \sigma(0)$ on 
$\D$ is again of the weights $(-1,\cdots,-1)$, 
we see that our algorithm in Theorem \eqref{bblcst:def} 
gives $\X'_{\rm min}$ which corresponds to 
$\varphi'_{\rm min}$. Therefore, by comparing with 
\cite[\S 3.3, Theorem 3]{dBS} and iterate the 
construction, we complete the proof. 

The case 
\eqref{cl3} is easy as the deformation space 
${\rm Def}(X=A_1\ni 0)$ is $1$-dimensional and smooth, 
which is realized as the $\C^*$-equivariant degeneration to 
the cone over quadric, which also serves as 
an example of 
(universal) negative weight deformation 
of Corollary \ref{Kur} 
in the next subsection. 
Each general fiber comes with the Stenzel metric 
\cite{Ste} for the complex quadric 
and 
it is easy to see that the construction in 
Theorem \ref{bblcst:def} gives each general fiber. 
Thus, by \cite[\S 4.2 Proposition 1]{dBS}, we complete the proof. 
\end{proof}

\subsubsection{Negative weight deformation as scaling down}\label{sec:nwd}

This subsection is not used in the main discussion and 
can be regarded as a complementary remark for those who 
know \cite{CH} or \cite{SunZhang}. 
Recall that the family $\mathcal{W}\to \A^1$ 
we constructed in Theorem \ref{bblcst:def} (see its 
\eqref{Kss.Song}) is an 
algebraic negative weight deformation (also called 
scale down degeneration) in the sense of \cite[\S 2]
{Od24b}. 
This small 
subsection aims to clarify the notion in the original 
more general setup over an affine toric variety. 
That is, we discuss the subtle 
relation between this algebraic notion in 
the general sense 
to 
the more original notion of negative weight deformation 
in the sense of \cite{CH}, which uses sequence of 
test configurations. 

This complements the easier direction assertion in 
\cite[Lemma 2.18(i)]{Od24b}, with which we 
partially answer the 
Open Question 1.4 (2) of \cite{CH}. 
Note that, in the context of the previous section, 
the degenerations 
$\X'_0\rightsquigarrow W$ and 
$\X''_0\rightsquigarrow C$ are both 
negative weight deformations by the construction. 

\begin{Thm}[Negative weight deformation]\label{negwtdeform}
\begin{enumerate}
\item (\cite[Lemma 2.18 (i)]{Od24b}) \label{1st}
Take an arbitrary 
algebraic negative weight deformation 
(cf., Definition 2.6 of \cite{Od24b})
of a Fano cone $T\curvearrowright(Y,\xi)$ as 
$p\colon \mathcal{Y}
\twoheadrightarrow U_{\tau}$, 
with $\xi\in \tau$. 
Then, the general fiber of $p$ 
is a negative $\xi$-weight deformation in the original sense of Definition 1.7 of \cite{CH}. 
\item \label{2nd}
Suppose that for a Fano cone $(T\curvearrowright Y,\xi)$, 
$Y$ has only an isolated singularity. 
Then, conversely to the above 
\eqref{1st}, 
any negative $\xi$-weight deformation $X$ of normal affine variety 
of $Y$ with $\xi$ in the sense of \cite{CH} 
arises in this manner. 

\end{enumerate}
\end{Thm}

These are also inspired by the examples of 
\cite{LiCY, CR, HN, Sze, Sze.uni, Chiu}. 

\begin{proof}
The above item \eqref{1st} is proved in 
\cite[Lemma 2.18(i)]{Od24b}, so it remains to prove 
\eqref{2nd}. 
We first confirm the 
existence of $T$-equivariant semi-universal 
deformation of $T\curvearrowright Y$ after Pinkham. 
The tangent space $\mathbb{T}^1$ of the 
deformation of $Y$ naturally admits 
$T$-action by 
\cite[Proposition 2.2]{Pink} and the arguments 
(2.5) (2.6) of {\it loc.cit}. The existence of the 
equivariant family follows from the arguments in the proof of Schlessinger 
\cite[Theorem 2.11]{Schlessinger}, as in \cite[(2.9), (2.10)]{Pink}. 
Indeed, the inductive construction of $J_q$ 
in (2.9) of {\it loc.cit} (or \cite[2.11]{Schlessinger}) 
extends without 
any change. We similarly take $F^0$ as a column vector of 
$T$-homogeneous equations $f_i (i=1,\cdots,N)$ of 
some $T$-equivariant embedding 
$Y\hookrightarrow \A^N$, 
and construct the lift $F^q$ whose entries are 
all $T$-equivariant by projection to 
$T$-eigenspace at the inductive step for each $q=0,1,2,\cdots$. 
Thus, the upshot of the necessary modification is 
simply to replace the degree 
(with respect to the $\G_m$-action) to the weight vector 
with respect to the $T$-action. 
Summarizing up, we obtained 
a $T$-equivariant semiuniversal deformation 
$\mathcal{U}\to {\rm Def}(Y)$ as a formal family. 
It can be retaken as an algebraic family by Artin's 
theorem \cite{Artin}, 
which we denote by the same symbol. 

We decompose $\mathbb{T}^1$ into its $T$-eigensubspaces
$\mathbb{T}^1_{\alpha}$, for $\alpha\in M:={\rm Hom}(T,\mathbb{G}_m)$, and put
\[
 \mathbb{T}^1_{\xi<0}:=
 \bigoplus_{\langle\xi,\alpha\rangle<0}\mathbb{T}^1_{\alpha}.
\]
Let ${\rm Def}^{-}_{\xi}(Y)\subset {\rm Def}(Y)$ be the corresponding
$\xi$-negative Kuranishi locus; equivalently, it is the attracting locus of
$[Y]$ for the contracting direction associated with $\xi$.  After shrinking
around $[Y]$, it is a finite-type affine $T$-scheme with a contracting
$T$-action and tangent space $\mathbb{T}^1_{\xi<0}$.  We set
\[
 \mathcal U^{-}_{\xi}:=
 \mathcal U\times_{{\rm Def}(Y)}{\rm Def}^{-}_{\xi}(Y).
\]

Now, we consider the negative $\xi$-deformation of $Y$ in the sense of \cite[Definition 1.7]{CH} 
and take the realizing sequence as $(X_i,\xi_i)_{i=1,2,\cdots}$. 
From the conditions (1) and (2) of  \cite[Definition 1.7]{CH} (also cf., its (4)), 
for sufficiently large $i$, an equivariant classifying map for the
deformation appearing in \cite[Definition 1.7]{CH} gives a point
$b_i\in {\rm Def}^{-}_{\xi}(Y)$ whose fiber is isomorphic to $X$.
The normalization of $\overline{T\cdot b_i}$ is an affine toric
variety $U_\tau$, with $\xi\in\tau$.  Pulling back
$\mathcal U^{-}_{\xi}$ to $U_\tau$ gives the desired algebraic negative
weight deformation. 
Thus we obtain the desired assertion of \eqref{2nd}. 

For later convenience, in the isolated weighted-homogeneous hypersurface
case $Y=V(f)$, choose homogeneous representatives $g_1,\ldots,g_s$ of the
$\xi$-negative part of the Jacobian ring of $f$.  The corresponding
semi-universal negative deformation is explicitly
\[
 ({\rm Def}_x(Y))_{\xi<0}=\A^s_{t_1,\ldots,t_s},\qquad
 \mathcal U^{-}_{\xi}=V\!\left(f+\sum_{i=1}^s t_i g_i\right).
\]
\end{proof}

Note that in the proof of above Theorem \ref{negwtdeform} 
\eqref{2nd}, 
the following is proved, which we extract for later 
convenience. 

\begin{cor}[Semi-universal negative weight deformation]\label{Kur}
For each Fano cone $(T\curvearrowright Y, \xi\in N_{\R})$ 
with the only isolated hypersurface singularity $y$, 
there is a semi-universal 
negative weight deformation 
$u_{\xi<0}\colon \mathcal{U}_{\xi <0}\twoheadrightarrow  
({\rm Def}_x(Y))_{\xi<0}\subset {\rm Def}_x(Y)$ 
over a finite type base with a good $T$-action 
$({\rm Def}_x(Y))_{\xi<0}$ (in the sense of 
\cite{LS}) 
such that any (algebraic) 
negative $\xi$-weight deformation 
can be obtained by pullbacks of $\pi_{\xi<0}$. 


\end{cor}

\vspace{4mm}

After the construction of the semi-universal negative deformation as above, 
it is natural to consider the corresponding moduli theory (cf., also \cite[\S 6]{SunZhang}). 
Here are some questions, which we hope to explore more in near future. 

\begin{defn}\label{moduli.McCY}
For the above Kuranishi space 
${\rm Def}_{\xi<0}(Y)$ of Corollary \ref{Kur} 
for the negative $\xi$-weight deformation in Corollary \ref{Kur}, 
note that there is a canonical $G:={\rm Aut}_T(Y)$-action 
(automorphism group of $Y$ as Fano cone) 
with canonical linearization which is trivial on the origin 
$0=[Y]$. 
If $({\rm Def}_{\xi<0}(Y)\setminus 0)\curvearrowleft G$ 
is prestable or semistable 
in the sense of \cite[Chapter I, \S 4]{GIT}, 
we define 
the corresponding GIT quotient 
$M^{\xi<0}_Y:=({\rm Def}_{\xi<0}(Y)\setminus 0)//G$. 
\end{defn}
Recall $G$ is reductive as proved in 
\cite[Appendix]{DSII} over $\C$, 
and for more general base field, it is algebraically proved 
\cite[Corollary 3.13]{Od24a}. 

\begin{Ques}
In the above setup of Definition \ref{moduli.McCY}, 
do the corresponding GIT quotients 
$$M^{\xi<0}_Y:=({\rm Def}_{\xi<0}(Y)\setminus 0)//G$$ 
admit moduli interpretation of affine $\partial \bar{\partial}$-exact  
complete Calabi-Yau metrics with 
Euclidean volume growth? 
Do the prestability or semistability condition 
always hold? 
\end{Ques}
\begin{Rem}
For instance, if $G=T=\G_m$, the latter is true and  $M_Y^{\xi<0}$ 
always comes with a canonical realization as a closed (projective) 
subscheme inside a weighted projective space. 
If $Y$ is toric with only isolated singularity, 
then it is also true by \cite{AltmannT1, AltmannInvent} 
(indeed, for $n\ge 4$, the tangent space is 
trivial). 
Note that the recent interesting non-existence result 
of asymptotically conical K\"ahler metric 
on $G_2$ type complex symmetric spaces $X$ of rank $2$ 
with a specified asymptotic (horospherical) cone $(Y,\xi)$ 
\cite[Corollary D]{Ngh} is due to 
the discrepancy of the (metric) Reeb vector field 
and the negative valuation i.e., 
$X\notin {\rm Def}_{\xi}^{-}(Y)$, 
hence is of different nature from what is 
discussed here. We thank T.T.Nghiem for 
the clarification on this point. 
\end{Rem}


\section{Alternative description via coordinates and weights}\label{sec:min1}


In this section, we discuss an alternative description of the 
minimal bubble in terms of coordinates and defining equations. 

We again work on $\pi\colon \X\to \Delta$ as 
the Setup \ref{Setup}. As in \cite[\S 2.2 Example 2.9]{Od24b}, 
we assume (without loss of generality) 
$x\in X$ and $0\in W$ are embedded into the affine space $0\in \A_{\mathbb{C}}^l$ and 
the isotrivial degeneration $X\rightsquigarrow W$ is obtained 
inside $\A^l_{\C[[t]]}$ 
by rescaling with weights $(w_1,\cdots,w_l)\in \R_{>0}^l$ for the 
coordinates $z_i$ of $\A_{\mathbb{C}}^l$. 
We also put a weight on $t$, the base coordinate, as $w_0\in 
\R_{>0}$. 


We further fix the extended embedding of $\X$ 
into $\Delta\times \A^l$ which maps the section $\sigma(\Delta)$ 
to the $0$-section $0\times \Delta(\simeq \Delta)$. 
Now, we consider $$\W:=\W(w_0;w_1,\cdots,w_l):=\Spec(gr_{w_0,w_1,\cdots,w_l}(\O_{\X,x})),$$ 
as a family over $\A^1_t$ where $t$ is the 
coordinate of $\Delta$, as in \cite{dBS}. 
Here, $gr_{w_0,w_1,\cdots,w_l}(\O_{\X,x})$ is the natural graded ring 
\begin{align}\label{gr.def}
\oplus_{\lambda\in \R_{\ge 0}}
\{\overline{f}\in \O_{\X,x}\mid f\in \O_{\Delta\times \A^l,0}, w(f)\ge \lambda\}/
\{\overline{f}\in \O_{\X,x}\mid f\in \O_{\Delta\times \A^l,0}, w(f)> \lambda \}, 
\end{align}
where for $f=\sum a_{i_0,i_1,\cdots,i_l}t^{i_0}\prod_{1\le j\le l}z_j^{i_j}$, the weight 
$w(f)$ of $f$ means $$\min \{w_0i_0+\sum_j w_j i_j \mid a_{i_0,i_1,\cdots,i_l}\neq 0\}$$
and $\overline{f}$ means the restriction. 
Note that this construction \eqref{gr.def} 
is a variant of the graded ring discussed in 
e.g., \cite{LXZ} but its finite typeness over 
$\mathbb{C}$ is straightforward in this setup. 
Indeed, if 
the ideal sheaf for $\X\subset \A^l \times \Delta$ 
is $I$, and $w(\mathcal{O}_{\Delta\times \A^l,0})$ 
consists of $\lambda_1<\lambda_2<\cdots$, if we set 
$\{f\in \O_{\Delta\times \A^l,0}\mid f\in \O_{\Delta\times \A^l,0}, w(f)\ge \lambda_i\}=:R_i$, 
we have 
\begin{align}
&\{\overline{f}\in \O_{\X,x}\mid f\in \O_{\Delta\times \A^l,0}, w(f)\ge \lambda_i\}/
\{\overline{f}\in \O_{\X,x}\mid f\in \O_{\Delta\times \A^l,0}, w(f)\ge \lambda_{i+1}\}\\ 
&=R_i/
(R_{i+1}
+
\{f\in \O_{\Delta\times \A^l,0}\mid f\in 
\O_{\Delta\times \A^l,0}, w(f)\ge \lambda_i \}\cap 
I)\twoheadleftarrow R_i/R_{i+1}, 
\end{align}
for each $i=1,2,\cdots$ 
and that 
${\rm gr}_{w_0,\cdots,w_l}(\mathcal{O}_{\A^l,0})
=\oplus_{i} R_i/R_{i+1}\simeq 
\mathbb{C}[z_1,\cdots,z_l,t]$ is a finite type algebra 
over $\mathbb{C}$. 

\begin{prop}[Critical weight]\label{bblcst:coord}
In the above context, 
\begin{enumerate}
\item \label{crit}
there is a natural
critical value $w_{\rm cri}\in \R_{\ge 0}$, such that 
\begin{itemize}
\item $w_0>w_{\rm cri}$, then $\W$  is canonically isomorphic to $\A^1_t\times W$ over $\A^1_t$. 
\item $w_0=w_{\rm cri}$, then $\W$  is an affine faithfully 
flat and $\G_m$-equivariant family with the central fiber $W$, where 
$\G_m$ acts on the base $\A_t^1$ with the weight $w_0$. 
\end{itemize}
\item \label{coord.val}
In particular, for $w_0\ge w_{\rm cri}$, the weights $w_0,w_1,\cdots,w_l$ induce 
a valuation on $\O_{\X,x}$. 
\end{enumerate}
\end{prop}
This proposition is inspired by \cite{Song, dBS} and the discussions with the authors. 
\begin{proof}
As we assumed in the beginning of this subsection, 
$\X_{\C[[t]]}\subset \A^l_{\C[[t]]}$ is cut out by some 
defining polynomials $F_1,\cdots,F_m$ of $z_1,\cdots,z_l,t$. We (again) take enough $F_i$'s to 
form a universal Gr\"obner basis of the generating ideal (\cite[\S 5.1]{BCRV}). 
Since $0\times \Delta\subset  \X$, 
these do not contain constant terms. 
Then $\W\subset \A_t^1\times \A^l$ is cut out by 
the initial terms ${\rm in}_{w_0,w_1,\cdots,w_l}(F_i)$ of $F_i$ with respect to the weights 
$w_0, w_1,\cdots,w_l$. Thus in particular, $\mathcal{W}$ 
admits the natural $\G_m$-action which equivariantly projects down to 
$\A_t^1$ with the weight $w_0$. 

We set $f_i:=F_i|_{t=0}$ and write 
\[
F_i=f_i+tg_i^{(1)}+t^2g_i^{(2)}+\cdots,
\]
with each $g_i^{(j)}\in \C[z_1,\cdots,z_l]$. 
Since $\X\to\Delta$ is flat, its defining ideal is
$t$-saturated; after dividing each $F_i$ by the maximal possible power of
$t$, we may and do assume $f_i\neq0$. 
Then, we define $s_i^j:=\dfrac{w(f_i)-w(g_i^{(j)})}{j} (j=1,2,\cdots)$ and 
\[w_{\rm cri}:=\max\{0,\max_{i\ge 1, j\ge 1, g_i^{(j)}\neq 0}\{s_i^j\}.\}\] 
Here, $w(-)$ denotes the weight of the polynomial of $z_1,\cdots,z_l$ 
with respect to $w_1,\cdots,w_l$. 
Then, from the construction, the initial term of $F_i$ is ${\rm in}_{w}(f_i)$ for any $i$ if $w_0>w_{\rm cri}$. 
For $w_0=w_{\rm cri}$, 
${\rm in}_{w_0,w_1,\cdots,w_l}(F_i)|_{t=0}={\rm in}_{w}(f_i)$ 
for each $i$ so that the central fiber of $\W$ is still $W$. 

For general weights $w_0,\ldots,w_l$, the associated graded ring need not
be flat over $\A^1$.  Here $\dim\mathcal W=n+1$ and $\dim W=n$, and, more
importantly, Theorem~\ref{kl}\eqref{kl2} below identifies the critical
family with the restriction of the affine faithfully flat family
$\widetilde{\mathcal X}_F$.  Hence no $t$-torsion occurs and
$\mathcal W(w_{\rm cri},w_1,\ldots,w_l)\to\A^1$ is faithfully flat.
Its generic fiber is integral, so its coordinate ring, which injects into
its localization at $t$, is integral.  Thus
$w_{\rm cri},w_1,\ldots,w_l$ defines a valuation on $\O_{\X,x}$
(cf.\ \cite[proof of Proposition 3.3(a) and comments]{Teissier}). 
This completes the proof of the item \eqref{coord.val}. 
\end{proof}

We now prepare a key observation 
to revisit Theorem \ref{bblcst:def} 
from an approach via coordinates and weights. 
Henceforth, we require more background from \cite{Od24b}. Indeed, 
note that a core of our proof of Theorem \ref{bblcst:def} 
is a general semistable reduction \cite[Theorem 1.1, 3.5]{Od24b} (compare our Figure \ref{fig:semistable-minimal-bubble} and \cite[Figure 2, right side]{Od24b}!). 
We give a quick review of its detail here; 
as our initial input, we had a parameter space 
$T\curvearrowright H$ and an arc inside it $\Delta\to H$,  
in the beginning of the proof of 
Theorem \ref{bblcst:def}, where the latter specifies 
$\X\to \Delta$ and its ($0$-)section $\sigma$. 

In our proof of Theorem \ref{bblcst:def}, 
as a result of applying \cite[Theorem 1.1, 3.5]{Od24b}, we implicitly 
obtained 
a regular rational polyhedral cone 
$\tau'\subset N \oplus \Z$, whose generators can be taken inside $N\times 
\R_{> 0}$ and the normal vector is 
$(-m,1)\in ((\tau^o)\cap M)\times \Z$. 
For them, {\it op.cit} Definition 3.1  introduces 
\[\widetilde{\rm ST}_{\C[[t]]}(\tau'):=\Spec \biggl(\dfrac{\C[[t]][\tau^{\vee}\cap M][X]}{(t-X\chi^{m})}\biggr),\]
where $X$ is a new variable and $\chi^m$ denotes the character of $T$ 
for $m\in M$. 
Then, as a consequence of applying {\it op.cit} 
Theorem 1.1, 3.5 in our proof of Theorem \ref{bblcst:def}, 
there is a $T$-equivariant morphism 
\[F\colon 
\widetilde{\rm ST}_{\mathbb{C}[[t]]}(\tau')\to H.\] 

We take a primitive integral vector
\[
 v'=(v,w_0)\in {\rm relint}(\tau')\cap(N\oplus\Z),
 \qquad w_0>0,
\]
which projects to $v\in {\rm relint}(\tau)\cap N$.  Writing $u_i\in M$
for the character of the semi-invariant coordinate $z_i$, we henceforth
write $w_i:=\langle u_i,v\rangle$; choosing $v$ in this cone does not
change the corresponding initial degeneration $W$. 

This $w_0$ is precisely the
critical value $w_{\rm cri}$ of Proposition~\ref{bblcst:coord}: the
central-fiber condition gives $w_0\ge w_{\rm cri}$, whereas
$w_0>w_{\rm cri}$ would give the trivial family $W\times\A^1$. 

We consider the rank $1$ algebraic subtorus of $T$ 
associated to $\Z v$ as $(\G_m\simeq )T_v \subset T$. 
Its analog for $v'=(v,w_0)$ is the following: 
we consider a ring homomorphism to a polynomial ring as 
\begin{align}\label{item:varphipull}
\dfrac{\C[\tau^\vee\cap M][X,t]}{(X\chi^m-t)}(\simeq \C[\tau^\vee\cap M][X])\to \C[s]
\end{align}
defined by $X\mapsto 1$, $\chi^{u}\mapsto s^{\langle u,v\rangle}$
for $u\in M$, and $t\mapsto s^{w_0}$. Since $(-m,1)$ is
orthogonal to $\tau'$, we have $\langle m,v\rangle=w_0$, so this map
respects the relation $X\chi^m=t$.  It naturally extends to
\[
\Gamma(\widetilde{\rm ST}_{\mathbb{C}[[t]]}(\tau'),
\mathcal O_{\widetilde{\rm ST}_{\mathbb{C}[[t]]}(\tau')})
\longrightarrow \C[[s]],
\]
and induces
\[
 \varphi_{v'}\colon {\rm Spec}\C[[s]]
 \longrightarrow \widetilde{\rm ST}_{\mathbb{C}[[t]]}(\tau').
\]
\begin{Thm}[Description via coordinates and weights]\label{kl}
In the above situation, we have the following: 
\begin{enumerate}
    \item \label{kl1} 
    Consider 
    the affine faithfully flat family 
    which corresponds to $F$ and denote it by 
     $\tilde{\X}_{F}\to \widetilde{\rm ST}_{\mathbb{C}[[t]]}(\tau')$. Then, we denote its base change to 
     ${\rm Spec}\C[[s]]$ by $\varphi_{v'}$ as 
     $$\X_{v'}\to {\rm Spec}\C[[s]].$$ 

Then, if $w_0=1$, $\X_{v'}$ is a Zariski open subscheme of 
the weighted blow up of $\X$ 
    for weights $(w_0=1,w_1,\cdots,w_l)$. This open chart is obtained by removing the strict transform of the central fiber. 
    For general $w_0\ge 1$, 
    $\X_{v'}$ is the normalization of the base change of the weighted 
    blow up of $\X$ for weights $(w_0,w_1,\cdots,w_l)$ 
    with respect to the $w_0$-th power map 
    \[\Spec \C[[s]]\xrightarrow{t\mapsto s^{w_0}} \Spec \C[[t]].\] 

    \item \label{kl2}
    The negative weight deformation (affine test configuration)  
$\mathcal{W}$ (cf., \eqref{Kss.Song} of Theorem \ref{bblcst:def}) 
    obtained as the     
    restriction of the family 
    $\tilde{\X}_{F}\to \widetilde{\rm ST}_{\mathbb{C}[[t]]}(\tau')$ 
    to the $1$-dimensional toric stratum $(\A^1\simeq )V(\tau')\subset \widetilde{\rm ST}_{\mathbb{C}[[t]]}(\tau')$, coincides with 
\[\mathcal{W}={\rm Spec}({\rm gr}_{(w_0,\cdots,w_l)}(\mathcal{O}_{\X,x}))\to 
\A^1.\] Hence, $\W$ is the 
affine cone of the exceptional divisor of weighted blow up, 
which appears in \eqref{kl1}. 
\end{enumerate}
\end{Thm}

The above theorem gives alternative description to 
$\X_{\rm min}'$ and $\mathcal{W}$ obtained in 
Theorem \ref{bblcst:def}. Our \eqref{kl1} 
answers a question of S.Donaldson to the author. 

\begin{proof}[proof of Theorem \ref{kl}]
Firstly, note that 
$\tilde{\X}_{F}\to \widetilde{\rm ST}_{\mathbb{C}[[t]]}(\tau')$ 
is affine and faithfully-flat finite type from the definition. 
We first prove \eqref{kl1}. 
Take the weights of $T_v$ on the coordinates
$z_1,\ldots,z_l$ of $\A^l$ to be $w_1,\ldots,w_l$. 
 Note that $\A_{\mathbb{C}[[t]]}^l$ has a regular system of parameters as 
    $t,z_1,\cdots,z_l$ so that one can discuss its weighted blow up 
    for any weights $(w_0,\cdots,w_l)$. 
    
If $w_0=1$, equivalently, 
$\X_{v'}$ is the restriction of the Zariski closure of 
$\X_{\C[[t]]}\setminus X$ in $$\A^l_{\frac{z_1}{t^{w_1}},\cdots,\frac{z_l}{t^{w_l}}}\times_{\Spec \C} \Spec \C[[t]].$$ Thus, if we denote the weighted blow up of $\X_{\C[[t]]}$ with weights 
$(w_0=1,w_1,\cdots,w_l)$ as 
$b\colon {\rm Bl}_{(w_0=1,w_1,\cdots,w_l)}\X_{\C[[t]]}\to 
\X_{\C[[t]]}$, we have 
$\X_{v'}={\rm Bl}_{(w_0=1,w_1,\cdots,w_l)}\X_{\C[[t]]}
\setminus (b^{-1}_* X)$. 
Hence the assertion of 
\eqref{kl1} holds in this case. 

For general $w_0$ case, 
let us use the $w_0$-th power map of the base 
    $\Spec \C[[s]] \xrightarrow{s\mapsto s^{w_0}=t} \Spec \C[[t]]$, 
for base change. 
    Then the normalization of its base change 
    ${\rm Bl}_{(w_0,\cdots,w_l)}(\A_{\mathbb{C}[[t]]}^l)\times_{\mathbb{C}[[t]]}\mathbb{C}[[s]]$ is ${\rm Bl}_{(1,w_1,\cdots,w_l)}(\A_{\mathbb{C}[[s]]}^l)$ as it is easy to see ring-theoretically. Here, the latter weighted blow up 
    means that with respect to the regular system of parameters 
    $s,z_1,\cdots,z_l$. 
    Therefore, as their closed subschemes, there is a finite morphism 
    \[\beta\colon {\rm Bl}_{(1,w_1,\cdots,w_l)}(\mathcal{X}_{\C[[t]]}\times_{\mathbb{C}[[t]]}\mathbb{C}[[s]]) 
    \to {\rm Bl}_{(w_0,w_1,\cdots,w_l)}(\mathcal{X}_{\C[[t]]})\times_{\mathbb{C}[[t]]}\mathbb{C}[[s]].\] 
This $\beta$ is the normalization.  Indeed, its source is flat over the
normal DVR $\C[[s]]$, and both its closed and generic fibers are
geometrically normal (the former by Theorem~\ref{bblcst:def}, and the
latter because the general fibers are klt).  Hence its source is normal;
moreover $\beta$ is finite and birational. 
    Hence, the assertion of \eqref{kl1} follows from the case of $w_0=1$ and we conclude its proof. 

Next, we prove \eqref{kl2}. 
The fact that $\mathcal{W}$, as the restriction of the 
family $\tilde{\X}_{F}\to \widetilde{\rm ST}_{\mathbb{C}[[t]]}(\tau')$ 
to $V(\tau')$, is a negative weight deformation 
follows from the existence of the natural $T$-action. 
Now we identify it with 
$\mathcal{W}\to \A^1$. Our intuition behind the following arguments is that 
for any $c\in \mathbb{C}$, admitting 
$c=0$, compare 
$\lim_{s\to 0}(s^{w_1},\ldots,s^{w_l})\cdot
X_{c s^{w_0}}$ and a fiber of $\mathcal W$ over $c$.  Here,
$X_{c s^{w_0}}$ denotes the fiber of $\pi$ over $c s^{w_0}$. 

To give precise arguments, 
we introduce a certain variant of 
$\varphi_{v'}$ as follows. 
Consider an extra variable $c$ and a ring homomorphism 
\begin{align}\label{item:varphiextend}
\Gamma(\widetilde{\rm ST}_{\mathbb{C}[[t]]}(\tau'),\mathcal{O}_{\widetilde{\rm ST}_{\mathbb{C}[[t]]}(\tau')})
=\dfrac{(\C[[t]])[\tau^\vee\cap M][X]}{(X\chi^m-t)}
\to (\C[c])[[s]], 
\end{align}
defined by 
$\chi^u\mapsto s^{\langle u,v\rangle}$ for $u\in M$, 
$X\mapsto c$, and $t\mapsto cs^{w_0}$.  Note that specialization 
$c\mapsto1$ gives \eqref{item:varphipull} after completion to 
$\C[[s]]$. 
Thus, \eqref{item:varphiextend} provides
\begin{align}\label{varphitilde}
\tilde{\varphi}\colon 
\Spec\bigl(\mathbb{C}[c][[s]]\bigr)\to \widetilde{\rm ST}_{\mathbb{C}[[t]]}(\tau').
\end{align}

We now prove \eqref{kl2}. Note that since $v\in{\rm relint}(\tau)$, 
the reduction of $\tilde{\varphi}$ modulo $s$ identifies 
$\Spec\C[c]$ with $V(\tau')$. 
As a pullback (the fiber product) of $\tilde{\X}_{F}\to 
{\widetilde{\rm ST}}_{\C[[t]]}(\tau')$ by $\tilde{\varphi}$, 
we obtain a family 
$\mathcal{Z}:=\tilde{\X}_{F}
\times_{\widetilde{\rm ST}_{\C[[t]]}(\tau')}
\Spec\bigl(\C[c][[s]]\bigr)
\to \Spec\bigl(\C[c][[s]]\bigr)$. 
Then, from the definition, $\mathcal{W}=\mathcal{Z}|_{s=0}$. 

We identify $\mathcal{W}$ by using $\mathcal{Z}$. 
After completing at $x$, which does not change the associated
graded ring since all the weights are positive, we can and do write
\[
 \widehat{\mathcal O}_{\mathcal X,x}
 =\C[[t,z_1,\ldots,z_l]]/I.
\]
Over $s\neq 0$, the family $\mathcal{Z}$ is obtained from
$\mathcal X_{t=cs^{w_0}}$ by change of variables $z_i$s to new 
$z'_i$s with the equations $z_i=s^{w_i}z'_i$ for $i=1,\cdots,l$. 
Hence its defining ideal over 
$P[s^{-1}]$ is 
\[
I(cs^{w_0},s^{w_1}z'_1,\ldots,s^{w_l}z'_l)P[s^{-1}], 
\]
where $P:=\C[c][[s]][z'_1,\ldots,z'_l].$ 
Here, $I(cs^{w_0},s^{w_1}z'_1,\ldots,s^{w_l}z'_l)$ denotes the ideal of $P$ 
defined by the substitution; 
$t\mapsto cs^{w_0}, z_i\mapsto 
s^{w_i}z'_i (i=1,\cdots,l)$. 
Thus the defining ideal of 
$\mathcal{Z}$ is 
\[
 \mathfrak J
 =
 \left(
 I(cs^{w_0},s^{w_1}z'_1,\ldots,s^{w_l}z'_l)P[s^{-1}]
 \right)\cap P,
\]
so that that of $\mathcal{W}$ is 
\[
 \mathfrak J|_{s=0}
 =
 \left(\left(
 I(cs^{w_0},s^{w_1}z'_1,\ldots,s^{w_l}z'_l)P[s^{-1}]
 \right)\cap P\right)|_{s=0},
\]
as an ideal of $\C[c,z_1,\cdots,z_l]$. 
As in the proof of Proposition \ref{bblcst:coord}, if we take the 
universal Gr\"obner basis, it is straightforward to see that 
$\mathfrak J|_{s=0}={\rm in}_{w_0,w_1,\cdots,w_l}(I)$. 
Thus, we obtain 
\begin{align}
\mathcal{W}=
 \Spec\bigl({\rm gr}_{(w_0,\cdots,w_l)}
 (\mathcal O_{\mathcal X,x})\bigr)
 \to \A^1_c, 
\end{align}
which completes the proof of \eqref{kl2}. 

\end{proof}

The above theorem is inspired by the 
explicit descriptions of 
the metric bubblings via weighted blow ups in some classical examples 
de Borbon-Spotti (\cite{dBS}). 
Note that in order to obtain the differential 
minimal bubble in general, we need a good choice of 
extension of the embedding $X$ into $\A^l$ (resp., $\X$ into 
$\Delta\times \A^l$)  
as the construction of \cite[Lemma 4.6]{Song} suggests. 



\subsubsection*{Acknowledgements}
The author heartily thanks 
S.Sun and C.Spotti for the helpful conversations over many 
years on this topic. He also thanks H.Blum, 
S.Boucksom, M.de Borbon, 
S.K.Donaldson, D.Halpern-Leistner, 
M.Jonsson, Y.Kusakabe, 
M.Mauri and T.T.Nghiem 
for helpful comments and discussions. 
During the research, the author was partially 
supported by Grant-in-Aid for Scientific Research (B) 23H01069, 21H00973, 
Grant-in-Aid for Scientific Research (A) 21H04429, 
20H00112, 
and Fund for the Promotion of Joint International Research (Fostering Joint International Research) 23KK0249. 
Part of this work is done during 
the author's stay at Paris, notably 
Institut de math\'ematiques de Jussieu - 
Paris Rive Gauche from March 2024 and 
Imperial College London in June 2024. 
The author received there above-mentioned beneficial 
comments, and nice pleasant working environments. 
We appreciate their kind hospitality. 

The original version of this paper appeared as arXiv:2406.14518. 
Following some suggestions of the referee, 
due to the intricate nature of the construction, we have made the 
exposition simpler by (e.g., assuming Gorenstein properties) 
eliminating some technical statements while keeping the core arguments. 
Further, we 
divided the original paper into two parts. This paper is the first part  
(and the second will be \cite{Odfuturebble}). Also, generative AI 
is used to help to find minor errors and typos 
which are all fixed and 
confirmed by the author. 



\vspace{5mm} \footnotesize \noindent
Contact: {\tt odaka@kurims.kyoto-u.ac.jp} \\
RIMS, Kyoto University, Kyoto 606-8285. JAPAN \\

\end{document}